\documentclass[a4paper]{article}

\usepackage{mypackages}
\usepackage{mycommands}
\usepackage{plot_options}

\usepackage[margin=1in]{geometry}
\usepackage{authblk}
\usepackage[square,numbers]{natbib}
\usepackage{graphicx}

\providecommand{\sep}{\textbullet\ }

\begin{document}


\title{Boundary Integral Formulation of the Cell-by-Cell Model of Cardiac Electrophysiology\footnote{This work was supported by the European High-Performance Computing Joint Undertaking EuroHPC under grant agreement No 955495 (MICROCARD) co-funded by the Horizon 2020 programme of the European Union (EU) and the Swiss State Secretariat for Education, Research and Innovation.}}

\author[1]{Giacomo Rosilho de Souza\thanks{giacomo.rosilhodesouza@usi.ch}}

\author[1,2]{Rolf Krause\thanks{rolf.krause@usi.ch}}

\author[3,1]{Simone Pezzuto\thanks{simone.pezzuto@unitn.it}}

\affil[1]{Center for Computational Medicine in Cardiology, Euler Institute, Università della Svizzera italiana, via G.~Buffi 13, 6900 Lugano, Switzerland}

\affil[2]{Faculty of Mathematics and Informatics, FernUni, Schinerstrasse 18, 3900, Brig, Switzerland}

\affil[3]{Department of Mathematics, Università di Trento, via Sommarive 14, 38123, Povo (Trento), Italy}

\maketitle

\begin{abstract}
We propose a boundary element method for the accurate solution of the cell-by-cell bidomain model of electrophysiology. The cell-by-cell model, also called Extracellular-Membrane-Intracellular (EMI) model, is a system of reaction-diffusion equations describing the evolution of the electric potential within each domain: intra- and extra-cellular space and the cellular membrane. The system is parabolic but degenerate because the time derivative is only in the membrane domain. In this work, we adopt a boundary-integral formulation for removing the degeneracy in the system and recast it to a parabolic equation on the membrane. The formulation is also numerically advantageous since the number of degrees of freedom is sensibly reduced compared to the original model. Specifically, we prove that the boundary-element discretization of the EMI model is equivalent to a system of ordinary differential equations, and we consider a time discretization based on the multirate explicit stabilized Runge--Kutta method. We numerically show that our scheme convergences exponentially in space for the single-cell case. We finally provide several numerical experiments of biological interest.
\end{abstract}



\section{Introduction}\label{sec:intro}
The human heart is composed of billions of electrically-active myocytes. Altogether, myocytes form a syncytium of cells that enables electrical and mechanical synchronization of the tissue~\cite{kleber2004}.
Cardiac myocytes are excitable cells that can react and transmit electric currents to communicate and coordinate their action. Electrical propagation depends on the conductive properties of the cytoplasm and the extracellular matrix. Cell-to-cell conduction occurs via gap junctions, permeable channel mostly distributed in the myocyte longitudinal direction. Myocyte excitability is due to hundreds of thousands ion channels embedded in the cellular membrane. The overall propagation of the cardiac action potential emerges from a balance of diffusion and transmembrane currents.

Mathematically, cardiac electrophysiology models are systems of reaction-diffusion equations. The reaction term results from transmembrane currents, which are voltage-dependent and regulated through a gating mechanism. Ion channel gating is typically modeled via Hodgkin--Huxley formalism, yielding a possibly large set of ordinary differential equations. The diffusion term captures the spatial coordination of the cardiac tissue. The state-of-the-art model is the bidomain system. Here, the intracellular and extracellular spaces are superimposed and homogenized~\cite{neu1993}. Patient-specific organ-scale simulations routinely employ the bidomain model and its monodomain approximation. Despite being physiologically accurate, the bidomain model fails to capture the sub-cellular tissue organization.
The cell-by-cell bidomain model accounts for the Extracellular-Membrane-Intracellular (EMI) tissue components as separated (yet coupled) entities~\cite{spach1995,TMR21}. The cell-by-cell model enables a more accurate description of tissue heterogeneities, a key aspect in heart failure and atrial fibrillation~\cite{schotten2011}.

The cell-by-cell model presents several challenges. First, its unusual mathematical formulation showing time dynamics at the boundaries, indeed it presents an ordinary differential equation (ODE) on the transmembrane boundary and a constraint on the gap junctions. Second, in addition to the natural stiffness introduced by the Laplacian, the ionic model introduces stiff nonlinear multiscale dynamics. Third, a full scale heart model would require billions of cells leading to an incredibly large system of equations.  Hence, advanced tailored methods must be designed to solve cell-by-cell models.

In the literature, cell-by-cell models have already been solved by means of the finite element or boundary element method. In the finite element community the problem has been tackled by Stinstra and collaborators~\cite{SHL09,SLH10,SPL07,SHL05} and more recently by Tveito and collaborators~\cite{JEG22,TJK17,TMR21,daversin2021} as well. In both cases the cell-by-cell model was employed to study the effects of the cells microscopic structure on macroscopic values as conductive velocity or effective tissue conductivity. Also, the cell-by-cell model was employed to derive the effective parameters for the bidomain model under different microstructural conditions~\cite{SLH10,spach1995,hand2009}. Bécue, Potse and Coudière \cite{bec16,BPC18,BPC17} compared different gap junctions modelizations and studied existence of solutions in \cite{becue}.
In the context of the boundary element method the model was solved only for very simple and structured geometries, for instance in \cite{FoS95,HLR92,LeR93} for a longitudinal array of non-touching cells, in \cite{VVV97} for two a two-cells model, and more recently, from a theoretical point of view, in \cite{HeJ18,HJA17} it was analyzed for the case of isolated cells.

In this paper we propose a spatial discretization of the cell-by-cell model based on the boundary element method (BEM) and reduce it to a single system of ODEs living only on the transmembrane boundary. The great advantage brought by the BEM is that only boundaries need to be discretized, leading to much smaller systems of equations compared to more traditional methods as finite elements or volumes. We stress that the approach presented here is easily adapted for different gap-junction boundary conditions~\cite{TMR21} or unbounded extracellular domains (the ``infinite bath'' approximation).  Also, any spatial discretization method for which Dirichlet-to-Neumann maps can be computed could be employed instead of the BEM. Compared to previous BEM approaches, our methodology is independent from the cells structure and reduces any problem to an ODE on the transmembrane boundary.

This paper is organized as follows. In \cref{sec:introex} we treat the simple case where only one myocyte cell is present, the purpose of this section is to introduce the needed tools and our approach in a simplified setting. In \cref{sec:fullprob} instead we discretize in space the full problem with an arbitrary number of cells, possibly in contact, and reduce it to a system of ODEs. Finally, in \cref{sec:num} we present some numerical results.

\section{The single-cell problem}\label{sec:introex}

The main purpose of this section is to introduce in a simplified setting the boundary integral formulation and the boundary element method (BEM) employed to discretize the full problem, done in \Cref{sec:fullprob}, and as well the approach used to reduce the space discrete problem into a system of ordinary differential equations (ODE).

\subsection{Problem formulation}

\begin{figure}
	\centering
		
	\tikzset{
		ring shading/.code args={from #1 at #2 to #3 at #4}{
			\def\colin{#1}
			\def\radin{#2}
			\def\colout{#3}
			\def\radout{#4}
			\pgfmathsetmacro{\proportion}{\radin/\radout}
			\pgfmathsetmacro{\outer}{.8818cm}
			\pgfmathsetmacro{\inner}{.8818cm*\proportion}
			\pgfmathsetmacro{\innerlow}{\inner-0.01pt}
			\pgfdeclareradialshading{ring}{\pgfpoint{0cm}{0cm}}%
			{
				color(0pt)=(white);
				color(\innerlow)=(white);
				color(\inner)=(#1);
				color(\outer)=(#3)
			}
			\pgfkeysalso{/tikz/shading=ring}
		},
	}

\begin{tikzpicture}	
	\fill[cellexterior,draw=black,thick] (0,0) ellipse (2.5 and 1.75);
	\fill[cellinterior,draw=cellmembrane,thick] (0,0) ellipse (1.5 and 0.75);

	\draw[thick,->] (1.5,0) -- (2,0) node[above,font=\large]  {$\bn_1$};
	\draw[thick,->] (1.5,0) -- (1,0) node[above,font=\large]  {$\bne$};
	\draw[thick,->] (-1.7678,    1.2374) -- (-2.0545,    1.6471) node[above,font=\large]  {$\bne$};

	\node[font=\large] at (0,0) {$\Omega_1$};
	\node[font=\large] at (-2,0) {$\Omega_0$};
	\node[font=\large] at (0.8,-0.95) {$\Gammam$};
	\node[font=\large] at (2.1,1.3) {$\Sigma$};
\end{tikzpicture}
	
	\caption{Geometrical setup of problem \cref{eq:cbcone}. The cell and the extra-cellular space are respectively denoted by $\Omega_1$ and $\Omega_0$. The cellular membrane is $\Gammam$.}
	\label{fig:geom}
\end{figure}
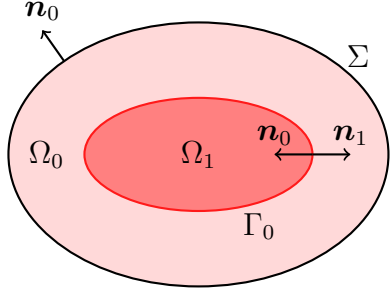

Here we consider the EMI model for a single cell, denoted by the bounded domain $\Omega_1\subset\Rd$ with $d=2$, embedded in the extracellular space, denoted by $\Omega_0 \subset\Rd$. See \cref{fig:geom} for a schematic representation of the single-cell problem. Specifically, we require that the intra- and extra-cellular domain do not overlap, that is $\Omega_0 \cap \Omega_1 = \emptyset$, and that they share a common boundary $\Gammam = \bar\Omega_0 \cap \bar\Omega_1$. The boundary $\Gammam$ represents the cellular membrane. The membrane model and temporal dynamic of the system, due to capacitative currents, is confined on $\Gammam$. We finally assume that $\Omega_0$ is bounded with exterior boundary $\Sigma = \partial\Omega_0 \setminus \Gammam$. Also, we define $\Gamma_1=\partial\Omega_1$ (note that for the single-cell problem $\Gamma_1=\Gamma_0$).
The single-cell problem reads as follows:
\begin{subequations}\label{eq:cbcone}
	\begin{empheq}[left={\empheqbiglbrace~}]{alignat=2}
		-\sigma_1 \Delta u_1 &= 0,
		&& \text{in $\Omega_1$}, \label{eq:cbcone_poisson_ui} \\
		-\sigmae \Delta \ue &= 0,
		&& \text{in $\Omegae$}, \label{eq:cbcone_poisson_ue} \\
		\sigma_1 \partial_{\bn_1} u_1 +\sigmae \partial_{\bne} \ue&=  0,
		&& \text{on $\Gammam$}, \label{eq:cbcone_cont_fluxes} \\
		\sigmae \partial_{\bne} \ue &=  \It(\Vm,z) = \Cm\partial_t\Vm + \Iion(\Vm,z) ,\qquad
		&& \text{on $\Gammam$},  \label{eq:cbcone_flux_ue_gammam}\\
		u_1 - \ue &=\Vm, && \text{on $\Gammam$},  \label{eq:cbcone_diff_u_eq_V}\\
		\partial_t z &= g(\Vm,z),
		&& \text{on $\Gammam$}, \label{eq:cbcone_ionic_var_eq}\\
            \sigmae \partial_{\bne} \ue &= 0 ,\qquad
		&& \text{on $\Sigma$},  \label{eq:cbcone_flux_ue_sigma}\\
	\end{empheq}
\end{subequations}
where $u_0(x,t)$, $u_1(x,t)$, and $\Vm(x,t)$ are respectively the intra-cellular, extra-cellular, and transmembrane electric potential, $\sigma_0$ and $\sigma_1$ are respectively the extra- and intra-cellular electric conductivity, $\Cm$ is the membrane capacitance, and $\bm{n}_i$, $i=0,1$ is the outwards normal. The ionic model is defined by $z(x,t)$, a vector of gating and concentration variables, its corresponding dynamic in \cref{eq:cbcone_ionic_var_eq}, and the ionic currents $\Iion(\Vm,z)$.

The global well-posedness of the problem~\eqref{eq:cbcone} in Sobolev spaces has been studied by \citet{matano2011}. The idea of the proof is similar to what we do here, in the sense that the authors recast \eqref{eq:cbcone} to an ODE on the interface $\Gammam$. The equation involves a pseudo-differential operator not dissimilar to the Dirichlet-Neumann map, as done below in the discrete settings with the operator $\mapV$. We are not aware of similar results for the multi-cell problem.

\subsection{Boundary integral formulation}\label{sec:cbem}

Problem in \cref{eq:cbcone} has already been tackled and carefully analyzed by \citet{HJA17}, where the BEM with a Galerkin approach was employed.
Here, we derive a boundary integral formulation of the unicellular problem \cref{eq:cbcone} in terms of trace operators and Poincaré--Steklov operators.
	
Let $\gammaut$ be the trace operator and $\gammaun$ the conormal derivative on the boundary $\Gamma_1$. More specifically, we introduce the operators as follows:
\begin{equation}\label{eq:defgammais}
	\begin{aligned}
		\gammaut &\colon H^1(\Omega_1) \to H^{1/2}(\Gamma_1), \quad &
		\gammaut u_1(\bx) &= \lim_{\Omega_1\ni \by\to \bx\in\Gamma_1} u_1(\by),  \\
		\gammaun &\colon H^1(\Omega_1) \to H^{-1/2}(\Gamma_1), &
		\gammaun u_1(\bx) &= \lim_{\Omega_1\ni \by\to \bx\in\Gamma_1}
		\langle \nabla u_1(\by), \bn_1 \rangle, 
	\end{aligned}
\end{equation}
where the limits on the right must hold for smooth enough $u_1$. 
Let $G(\bx,\by)$ be the fundamental solution of the Laplacian in $\Rd$, the Green representation formula for $u_1$ satisfying \cref{eq:cbcone_poisson_ui} implies that
\begin{equation}\label{eq:bir}
	\begin{aligned}
		u_1(\bx) &= \int_{\Gamma_1} \gamma^1_{t,\by} G(\bx,\by) \gammaun u_1(\by) \dif s_{\by}
		- \int_{\Gamma_1} \gamma^1_{n,\by} G(\bx,\by) \gammaut u_1(\by) \dif s_{\by}, && \bx\in\Omega_1
	\end{aligned}
\end{equation}
where $\by$ in $\gamma^1_{t,\by},\gamma^1_{n,\by}$ means that the operators are applied to the second variable of $G(\bx,\by)$. Taking the trace $\gammaut$ of \cref{eq:bir} we obtain
\begin{equation}\label{eq:traceGreen}
		\gammaut u_1 = \VV_1 \gammaun u_1-(\KK_1 -\tfrac{1}{2} I) \gammaut u_1,
\end{equation}
where $I$ is the identity operator, $\VV_1$ and $\KK_1$ are the single and double layer operators defined by
\begin{equation}\label{eq:sdop}
	\begin{aligned}
		\VV_1 &\colon H^{-1/2}(\Gamma_1) \to H^{1/2}(\Gamma_1), \quad &
		\VV_1 \rho(\bx) &= \int_{\Gamma_1} \gamma^1_{t,\by}G(\bx,\by)\rho(\by)\dif s_{\by}, \quad & 
		\bx & \in\Gamma_1, \\
		\KK_1 &\colon H^{1/2}(\Gamma_1) \to H^{1/2}(\Gamma_1), &
		\KK_1 \rho(\bx) &= \int_{\Gamma_1} \gamma^1_{n,\by} G(\bx,\by) \rho(\by)\dif s_{\by}, & 
		\bx & \in\Gamma_1, 
	\end{aligned}
\end{equation}
and the term $\tfrac{1}{2} I$ in \cref{eq:traceGreen} comes from the jump of the double layer potential as $\Omega_1\ni\bx\to\by\in\Gamma_1$.
Problem \cref{eq:traceGreen} can be rewritten as
\begin{equation}\label{eq:preDtN}
	\VV_1 \gammaun u_1=(\KK_1 +\tfrac{1}{2} I) \gammaut u_1,
\end{equation}
or employing the Poincaré--Steklov operator (Dirichlet-to-Neumann map):
\begin{equation}\label{eq:DtN}
\PP_1:H^{1/2}(\Gamma_1)\rightarrow H^{-1/2}(\Gamma_1),\qquad \PP_1\coloneqq\VV_1^{-1}(\KK_1 +\tfrac{1}{2} I),
\end{equation}
then \cref{eq:preDtN} becomes
\begin{equation}\label{eq:contPSi}
	\gammaun u_1=\PP_1 \gammaut u_1.
\end{equation}
The Poincaré--Steklov operator $\PP_1$ in \cref{eq:DtN} is known to be symmetric \cite[Section 3.7]{SaS11}.
Now, let $\gammaet u_0$, $\gammaen u_0$ be the trace and conormal derivative of $u_0$ on $\Gamma_0=\Gamma_1$, respectively.
In order to derive a Dirichlet-to-Neumann map $\PP_0:H^{1/2}(\Gamma_0)\rightarrow H^{-1/2}(\Gamma_0)$ in $\Omega_0$, hence
\begin{equation}\label{eq:contPSe}
	\gammaen u_0=\PP_0 \gammaet u_0,
\end{equation}
we need to take into account the boundary condition \eqref{eq:cbcone_flux_ue_sigma} on the external boundary $\Sigma$ of $\Omega_0$. In order to alleviate the presentation we postpone the derivation of $\PP_0$ to \ref{app:PSe}. However, we would like to note that $\PP_0$ remains symmetric.

Due to the Green representation formula, \cref{eq:cbcone_poisson_ui,eq:cbcone_poisson_ue} can be dropped from \cref{eq:cbcone}. Also, \cref{eq:cbcone_flux_ue_sigma} is encoded into the definition of $\PP_0$ (see \ref{app:PSe}). Finally, the boundary integral formulation of \cref{eq:cbcone} is
\begin{subequations}\label{eq:cbconebndint}
	\begin{empheq}[left={\empheqbiglbrace~}]{alignat=1}
		\sigma_1 \PP_1\gammaut u_1 +\sigmae \PP_0\gammaet u_0&= 0,  \label{eq:cbconebndinta} \\
		\gammaut u_1 - \gammaet u_0 &= V_0,  \label{eq:cbconebndintc}  \\
		\sigmae\PP_0 \gammaet u_0 &=  \It(V_0,z), \label{eq:cbconebndintb} \\
		\partial_t z &= g(V_0,z). \label{eq:cbconebndintd}
	\end{empheq}
\end{subequations}

\subsection{Spatial discretization of the unicellular problem}
We adopt the collocation BEM as spatial discretization scheme. Boundary element methods have less degrees of freedom than other standard techniques, while the collocation approach yields lower dimensional boundary integrals than the variational method and hence faster computations. For extensive presentations on the BEM we refer to \cite{Kre89,SaS11,Ste08}.

We place $M$ collocation points $\bx_j$, $j=1,\ldots,M$, on $\Gamma_1$ in a counterclockwise order. Then we compute a smooth parametrization $\gamma_{\Gamma_1}:[0,1)\rightarrow \Gamma_1$ satisfying 
\begin{equation}\label{eq:interpgamma}
	\gamma_{\Gamma_1}(t_j)=\bx_j\quad j=1,\ldots,M,
\end{equation}
where $\{t_j\}_{j=1}^M\subset[0,1)$ is an increasing sequence. (The parametrization $\gamma_{\Gamma_1}$ is computed with Fourier interpolation. For the unicellular problem, we could define $\gamma_{\Gamma_1}(t)$ first and then set $\bx_j$ as in \cref{eq:interpgamma}. However, this is not possible for the multi cell problems.) Finally, we represent $\gammaut u_1$, $\gammaun u_1$ as
\begin{equation}\label{eq:discgammasui}
	\begin{aligned}
		\gammaut u_1 (\gamma_{\Gamma_1}(t))&= \sum_{j=1}^{M} u_1^j L_j(t), & \gammaun u_1 (\gamma_{\Gamma_1}(t))&= \sum_{j=1}^{M} \tilde{u}_1^j L_j(t),
	\end{aligned}
\end{equation}
where $L_j(t)$ are trigonometric Lagrange polynomials satisfying $L_j(t_k)=\delta_{jk}$ for $j,k=1,\ldots,M$. Instead of \cref{eq:preDtN} we solve the weaker form
\begin{equation}
	\VV_1 \gammaun u_1(\bx_k)=(\KK_1 +\tfrac{1}{2} I) \gammaut u_1(\bx_k), \qquad k=1,\ldots,M,
\end{equation}
with $\bx_k=\gamma_{\Gamma_1}(t_k)$, which is equivalent to
\begin{equation}\label{eq:discDtN}
	\sum_{j=1}^{M} \tilde{u}_1^j \VV_1 L_j(\gamma_{\Gamma_1}^{-1}(\bx_k)) =\sum_{j=1}^{M} u_1^j (\KK_1 +\tfrac{1}{2} I) L_j(\gamma_{\Gamma_1}^{-1}(\bx_k)), \qquad k=1,\ldots,M,
\end{equation}
and hence the linear system
\begin{equation}\label{eq:prePSi}
	V_1 \btuu = (K_1+\tfrac{1}{2}I)\bui[1],
\end{equation}
with $\bui[1],\btuu$ the vectors of coefficients $u_1^j,\tilde{u}_1^j$, respectively, and 
\begin{align*}
	(V_1)_{kj} &\coloneqq \VV_1(L_j\circ \gamma_{\Gamma_1}^{-1})(\bx_k)= \int_{\Gamma_1} \gamma^1_{t,\by}G(\bx_k,\by) L_j(\gamma_{\Gamma_1}^{-1}(\by))\dif s_{\by}
	= \int_0^1 \gamma^1_{t,\by}G(\bx_k,\gamma_{\Gamma_1}(t)) L_j(t) \Vert \gamma_{\Gamma_1}'(t)\Vert \dif t,\\
	(K_1)_{kj} &\coloneqq \KK_1 (L_j\circ\gamma_{\Gamma_1}^{-1})(\bx_k)= \int_{\Gamma_1} \gamma^1_{n,\by} G(\bx_k,\by) L_j(\gamma_{\Gamma_1}^{-1}(\by))\dif s_{\by} 
	= \int_0^1 \gamma^1_{n,\by} G(\bx_k,\gamma_{\Gamma_1}(t)) L_j(t)  \Vert \gamma_{\Gamma_1}'(t)\Vert \dif t.
\end{align*}
The matrix coefficients $(V_1)_{kj}$, $(K_1)_{kj}$ must be computed with special care due to the singularities in the fundamental solution $G(\bx,\by)$ and its derivatives as $\gamma_{\Gamma_1}(t)\to \bx_k$, we refer to \cite{GKM21,Kre89} for the details.

Note that $\bui[1]$ and $\btuu$ are the vectors of coordinates of $\gammaut u_1$ and $\gammaun u_1$, respectively, and that from \cref{eq:prePSi} follow the discrete version of \cref{eq:contPSi}
\begin{equation}\label{eq:PSi}
	\begin{aligned}
		\btuu &=\PSi[1] \bui[1], & \PSi[1] & \coloneqq (V_1)^{-1}(K_1+\tfrac{1}{2}I),
	\end{aligned}
\end{equation}
where $\PSi[1]$ is the discrete Poincaré--Steklov operator (Dirichlet-to-Neumann map) in $\Omega_1$. Similarly, in \ref{app:PSe} we derive the discrete version of \cref{eq:contPSe} and obtain
\begin{equation}\label{eq:PSe}
	\begin{aligned}
		\btuz &=\PSi[0] \bui[0],
	\end{aligned}
\end{equation}
with $\bue$ and $\btue$ the vectors of coordinates of $\gammaet\ue$ and $\gammaen\ue$, respectively.

Finally, the space discretization of boundary integral formulation \cref{eq:cbconebndint} is
\begin{subequations}\label{eq:cbconedisc}
	\begin{empheq}[left={\empheqbiglbrace~}]{alignat=1}
		\sigma_1 \PSi[1]\bui[1] +\sigmae \PSe\bue&= 0,  \label{eq:cbconedisca} \\
		\bui[1] - \bue &= \bVm,  \label{eq:cbconediscc}  \\
			\sigmae\PSe \bue &=  \It(\bVm,\bz), \label{eq:cbconediscb} \\
		\bz' &= g(\bVm,\bz), \label{eq:cbconediscd}
	\end{empheq}
\end{subequations}
where $\bVm(t)\in\Rb^M$ is the vector whose coefficients represent $\Vm(\bx_j,t)$ and analogously for $\bz$. The right-hand sides $\It$ and $g$ are applied to $\bVm,\bz$ component wise. If needed, the solution $u_1$ satisfying \cref{eq:cbcone_poisson_ui} is approximated via the Green identity \cref{eq:bir} and  \cref{eq:discgammasui}, \cref{eq:PSi}. We proceed similarly for $\ue$.

\subsection{The Lagrange multipliers approach for the unicellular problem}\label{sec:ex_lag}
Now we solve \cref{eq:cbconedisc} and to do so we employ the Lagrange multiplier method. 
We adopted this technique mainly for pedagogical reasons in regard of what will be presented in \cref{sec:fullprob}, since for the unicellular problem \cref{eq:cbcone} a more direct	 approach could be used.

In the remaining of this section we construct the linear map 
\begin{equation}\label{eq:varphi}
	\begin{aligned}
		\mapV &:\Rb^M\rightarrow\Rb^M,\qquad & \mapV(\bVm)=\sigmae\PSe\bue,
	\end{aligned}
\end{equation}
where $\bue$ satisfies \cref{eq:cbconedisca,eq:cbconediscc} (behind the scenes $\bui[1]$ is computed as well, but it is not needed as output of $\mapV$). Inserting \cref{eq:varphi} and $\It(\bVm,\bz)=\Cm\bVm' + \Iion(\bVm,\bz)$ in \cref{eq:cbconediscb,eq:cbconediscd} the problem reduces to the ODE
\begin{subequations}\label{eq:cbconeode}
	\begin{empheq}[left={\empheqbiglbrace~}]{alignat=1}
		\Cm\bVm' + \Iion(\bVm,\bz)&=\mapV(\bVm) ,\\
		\bz' &= g(\bVm,\bz),
	\end{empheq}
\end{subequations}
which can be integrated by any suitable time marching scheme. \cref{eq:cbconeode} has the same structure as the one derived in \cite{HJA17}, where a BEM for the cell-by-cell model without gap junctions is derived.

\Cref{thm:unicell} is the unicellular version of the more general \cref{thm:phii} below, which in turn takes inspiration from the work in \cite{LaS03}. We also remark that the Theorem is independent on the spatial discretization. For instance, a finite element discretization may be recast to \cref{eq:cbconedisc} by static condensation, that it by explicitly computing the discrete Poincar\'e--Steklov operator.

\begin{theorem}\label{thm:unicell}
	The linear map $\mapV$ from \cref{eq:varphi} satisfies $\mapV(\bVm)=\blambda$, with $\blambda\in\Rb^M$ and $\beta_1\in\Rb$ solution to 
	\begin{equation}\label{eq:finallinsysone}
		\begin{pmatrix}
			F & G \\ G^\top & 0 
		\end{pmatrix}
		\begin{pmatrix}
			\blambda \\ \beta_1
		\end{pmatrix}
		=
		\begin{pmatrix}
			\bVm \\ 0
		\end{pmatrix}.
	\end{equation}
The matrices $F\in \Rb^{M\times M}$, $G\in \Rb^M$ are defined by
\begin{equation}\label{eq:unicellFG}
F=-(\sigma_1^{-1}(\PSi[1]^+)^{-1}+\sigmae^{-1}(\PSe^+)^{-1}), \qquad G=\be,
\end{equation}
$\be\in\Rb^M$ is the vector of ones and 
\begin{equation}
		\PSe^+ = \PSe+\alpha_0 \be\be^\top, \qquad 	\PSi[1]^+ = \PSi[1]+\alpha_1 \be\be^\top,
\end{equation}
with $\alpha_0,\alpha_1>0$. If needed, $\bui[1],\bue$ are computed with
\begin{equation}
	\bue = \sigmae^{-1}(\PSe^+)^{-1}\blambda, \qquad 
	\bui[1] = -\sigma_1^{-1}(\PSi[1]^+)^{-1}\blambda + \beta_1\be.
\end{equation}
\end{theorem}
Hence, when solving the ODE system \cref{eq:cbconeode} with a time integration scheme, every time that $\mapV(\bVm)$ needs to be evaluated system \cref{eq:finallinsysone} is solved and $\mapV(\bVm)=\blambda$ is inserted in \cref{eq:cbconeode}. 
\begin{proof}[Proof of \cref{thm:unicell}.]
Vectors $\bui[1]$, $\bue$ are solutions to the smaller system in \cref{eq:cbconedisca,eq:cbconediscc}.
Since $\PSe$, $\PSi[1]$ are symmetric, $\bui[1]$, $\bue$ are also solution to the constrained minimization problem
\begin{equation}\label{eq:minprob}
			\min_{\bui[1],\bue} \frac{\sigmae}{2}\langle\PSe\bue,\bue\rangle +\frac{\sigma_1}{2}\langle \PSi[1]\bui[1],\bui[1]\rangle
			\quad  \text{with} \quad
			\bui[1]-\bue=\bVm,
\end{equation}
where $\langle\cdot,\cdot\rangle$ is the Euclidean inner product in $\Rb^M$. Let 
\begin{equation}\label{eq:lagrangian}
	\mathcal{L}(\bue,\bui[1],\blambda) = \frac{\sigmae}{2}\langle\PSe\bue,\bue\rangle + \frac{\sigma_1}{2}\langle \PSi[1]\bui[1],\bui[1]\rangle +\langle\bui[1]-\bue-\bVm,\blambda\rangle
\end{equation}
be the Lagrangian function, imposing $\nabla\mathcal{L}(\bue,\bui[1],\blambda)=\bm{0}$ yields
\begin{equation}\label{eq:linsysone}
		\sigmae\PSe \bue - \blambda=\bm{0},\qquad
		\sigma_1\PSi[1] \bui[1] + \blambda= \bm{0},\qquad
		\bui[1]-\bue =\bVm.
\end{equation}
Note that $\PSe \be=\bm{0}$, with $\be\in\Rb^M$ a vector of ones, thus $\PSe$ is singular and the first equation of \eqref{eq:linsysone} $\sigmae\PSe\bue=\blambda$ has a solution only if
\begin{equation}\label{eq:condlambda}
	\langle\blambda,\be\rangle=0.
\end{equation}
Let 
\begin{equation}\label{eq:defPSep}
	\PSe^+ = \PSe+\alpha_0 \be\be^\top
\end{equation}
with $\alpha_0>0$, then $\PSe^+$ is invertible and it can be verified that if \cref{eq:condlambda} holds then
\begin{equation}
	\bue = \sigmae^{-1}(\PSe^+)^{-1}\blambda + \beta_0\be, \qquad \beta_0\in\Rb,
\end{equation}
is solution to $\sigmae\PSe\bue=\blambda$. Similarly, 
\begin{equation}
	\bui[1] = -\sigma_1^{-1}(\PSi[1]^+)^{-1}\blambda + \beta_1\be, \qquad \beta_1\in\Rb,
\end{equation}
is solution to $\sigma_1\PSi[1] \bui[1] =- \blambda$, with $\PSi[1]^+$ defined analogously to $\PSe^+$. Note that if $\bue,\bui[1]$ are solutions to \cref{eq:cbconedisca,eq:cbconediscc} then $\bue+C\be$, $\bui[1]+C\be$ are solutions for all $C\in\Rb$. We choose to fix such free constant by setting $\beta_0=0$, which implies $\sum_{j=1}^M\ue^j=0$. The last equation of \cref{eq:linsysone} yields
\begin{equation}
	\begin{aligned}
		\bVm &= \bui[1]-\bue 
		=-(\sigma_1^{-1}(\PSi[1]^+)^{-1}+\sigmae^{-1}(\PSe^+)^{-1})\blambda + \beta_1\be 
		= F\blambda + G\beta_1,
	\end{aligned}
\end{equation}
with matrices $F$, $G$ as in \cref{eq:unicellFG}. Together with \cref{eq:condlambda} it yields system \cref{eq:finallinsysone}.
\end{proof}

\begin{remark}\label{rem:unbouned}
	Note that the content of this section is readily adapted to a problem \cref{eq:cbcone} with unbounded domain $\Omegae$, hence without boundary condition \eqref{eq:cbcone_flux_ue_sigma}. In that case, $\PSe$ is derived analogously to $\PSi[1]$. However, $\PSe$ would be non singular hence in \cref{thm:unicell} we consider \cref{eq:defPSep} with $\alpha_0=0$. Condition \cref{eq:condlambda} is still required for the existence of a solution to $\sigma_1\PSi[1]\bui[1]=-\blambda$.
\end{remark}

\section{Discretization of the full cell-by-cell model}\label{sec:fullprob}

\begin{figure}
	\centering
		
	\tikzset{
		ring shading/.code args={from #1 at #2 to #3 at #4}{
			\def\colin{#1}
			\def\radin{#2}
			\def\colout{#3}
			\def\radout{#4}
			\pgfmathsetmacro{\proportion}{\radin/\radout}
			\pgfmathsetmacro{\outer}{.8818cm}
			\pgfmathsetmacro{\inner}{.8818cm*\proportion}
			\pgfmathsetmacro{\innerlow}{\inner-0.01pt}
			\pgfdeclareradialshading{ring}{\pgfpoint{0cm}{0cm}}%
			{
				color(0pt)=(white);
				color(\innerlow)=(white);
				color(\inner)=(#1);
				color(\outer)=(#3)
			}
			\pgfkeysalso{/tikz/shading=ring}
		},
	}

\begin{tikzpicture}
	\fill[cellexterior,draw=black,thick] (5.5,1.1) ellipse (7.5 and 2.75);
	
	\draw [rounded corners,fill=cellinterior,draw=black] (-0.5,1.2)--(1.5,1.2)--(2,2.2)--(0,2.2)--cycle;
	\draw [rounded corners,fill=cellinterior,draw=cellmembrane] (0.5,0)--(2,0)--(2.5,1)--(1,1)--cycle;
	\draw [rounded corners,fill=cellinterior,draw=cellmembrane] (2,0)--(4.5,-0.5)--(5,0.5)--(4,0.5)--(5.5,1)--(6,2)--(4.5,2)--(2.5,1)--cycle;
	\draw [rounded corners,fill=cellinterior,draw=cellmembrane] (5.5,1)--(8,1)--(8.5,0.5)--(10,0.5)--(9.5,1)--(11,1)--(11.5,2)--(8.5,2)--(6,2)--cycle;
	\draw [rounded corners,fill=cellinterior,draw=cellmembrane] (4.5,-0.5)--(7,-0.5)--(7.5,0.5)--(5,0.5)--cycle;
	\draw [rounded corners,fill=cellinterior,draw=cellmembrane] (7,-0.5)--(10,-0.5)--(10.5,0.5)--(7.5,0.5)--cycle;
	
	\draw[very thick,colorgammaij] (5.55,1.1)--(5.95,1.9);
	\node[colorgammaij,font=\large] at (6.2,1.5) {$\Gamma_{43}$};
	
	\draw[very thick,colornormal,->] (9.5,2)--(9.5,2.5)  node[right,font=\large]  {$\bn_4$};
	\draw[very thick,colornormal,->] (3.5,1.5)--(3.75,1.0)  node[above right,font=\large]  {$\bne$};
	\draw[very thick,colornormal,->] (10.8033 ,-0.8445)--(10.9754,   -1.3140)  node[right,font=\large]  {$\bne$};
	
	\node[font=\large,black] at (5.5,3) {$\Omegae$};
	\node[font=\large,black] at (0.75,1.7) {$\Omega_1$};
	\node[font=\large,black] at (1.5,0.5) {$\Omega_2$};
	\node[font=\large,black] at (3,0.5) {$\Omega_3$};
	\node[font=\large,black] at (8,1.5) {$\Omega_4$};
	\node[font=\large,black] at (6,0) {$\Omega_5$};
	\node[font=\large,black] at (8.75,0) {$\Omega_6$};
	\node[font=\large,black] at (1.3,2.4) {$\Gamma_1$};
	\node[font=\large,black] at (-0.3,-1) {$\Sigma$};
\end{tikzpicture}
	
	\caption{Illustration of problem \cref{eq:cbc}.}
	\label{fig:cbc}
\end{figure}
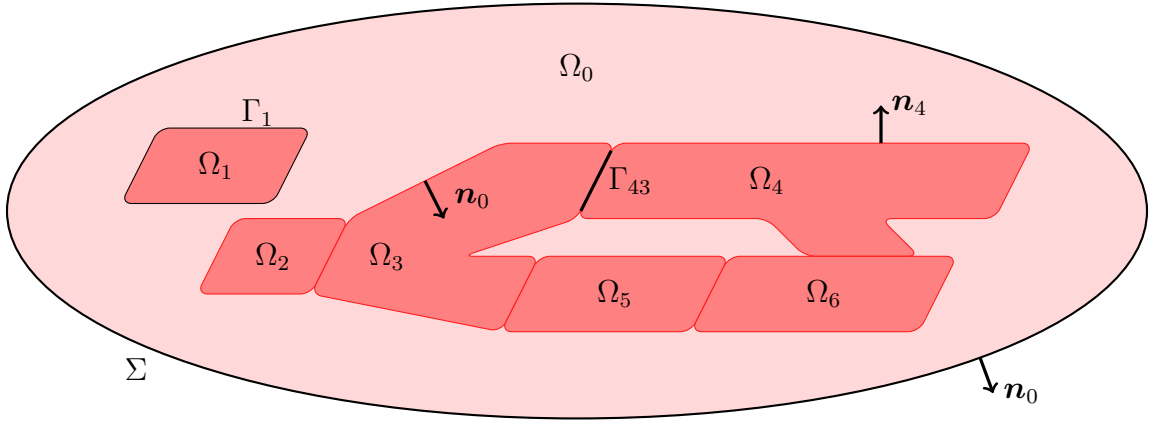

We introduce here the general cell-by-cell model. We consider an extracellular domain $\Omega_0\subset\mathbb{R}^d$, $d\ge 2$, an intracellular domain $\OmegaI\subset\mathbb{R}^d$, and an interface domain $\Gamma_0 = \bar{\Omega}_0\cap\bar{\Omega}_\mathrm{I}$. (See Figure~\ref{fig:cbc} for a graphical illustration of the model.) We suppose that $\OmegaI$ and $\Omega_0$ are disjoint and we denote by $\Omega$ the whole tissue, $\Omega = \OmegaI \cup \Omega_0 \cup \Gamma_0$. The domain $\Omega$ is always assumed connected and bounded, with $\Sigma = \partial\Omega$. For the sake of simplicity, $\partial\Omega_0 \setminus \Gamma_0 = \Sigma$ and $\partial\OmegaI = \Gamma_0$, that is the exterior boundary of $\Omega$ always corresponds to the extracellular matrix. Next, the intracellular space is described by the union of disjoint cells, denoted by $\Omega_i$, $i = 1,\ldots,N$. (Conveniently, $\Omega_i$ for $i=0$ corresponds to the extracellular space.) Thus, $\OmegaI = \bigcup_{i=1}^N \Omega_i$. We denote $\Gamma_i=\partial\Omega_i$, $i=1,\ldots,N$.
The cell-to-cell interconnections are denoted by $\Gamma_{ij} = \Gamma_i\cap\Gamma_j$, $0\leq j<i\leq N$. Note that the boundary of each cell is either in contact with another cell or with the extracellular space.
The cell-by-cell model reads as follows:
\begin{subequations}\label{eq:cbc}
	\begin{empheq}[left={\empheqbiglbrace~}]{alignat=2}
		-\sigmai \Delta \ui &= 0,
		&& \text{in $\Omegai$, $i=0,\ldots,N$}, \label{eq:cbc_poisson_ui} \\
		u_i - u_0 &=V_{i0}, 
		&& \text{on $\Gamma_{i0}$ for $1\leq i\leq N$},   \label{eq:cbc_diff_u_eq_V}\\
		-\sigmai \partial_{\bni} \ui &=  \Cm\partial_t V_{i0} + \Iion(V_{i0},z_i),
		&& \text{on $\Gamma_{i0}$ for $1\leq i\leq N$},   \label{eq:cbc_flux_ui_gammaim} \\
		-\sigmae \partial_{\bne} \ue &= - \Cm\partial_t V_{i0} - \Iion(V_{i0},z_i) ,\qquad
		&& \text{on $\Gamma_{i0}$ for $1\leq i\leq N$},   \label{eq:cbc_flux_ue_gammam} \\
		\partial_t z_i &= g(V_{i0},z_i),
		&& \text{on $\Gammam$}, \label{eq:cbc_ionic_var_eq} \\
		-\sigmai \partial_{\bni} \ui &=  \kappa (\ui-\uj),
		&& \text{on $\Gammaij$ for $1\leq j<i\leq N$},  \label{eq:cbc_permeability}\\
            -\sigmae \partial_{\bne} \ue &= 0 ,\qquad
		&& \text{on $\Sigma$}.   \label{eq:cbc_flux_ue_sigma}\\
	\end{empheq}
\end{subequations}
The constant conductivities are $\sigmai$, $i=0,\ldots,N$. The gap junctions (intercellular connections) are represented by $\Gammaij$ for $1\leq j<i\leq N$, with permeability $\kappa$. The normals $\bni$ point outwards to $\Omegai$. 
The intracellular potentials are $\ui$ for $i=1,\ldots,N$, the extracellular potential is $\ue$ and $V_{i0}$ is the transmembrane potential on $\Gamma_{i0}$. The membrane electric capacitance is $\Cm$ and $\Iion$ represents the sum of ionic currents. The transmembrane potential $V_{i0}$ is regulated by the ionic currents, which in turn depend on ionic concentrations and their transmembrane fluxes through ion channels, which are governed by gating variables. Ion concentrations and gating variables are represented by $z_i$ and the pair $\Iion,g$ describe the membrane ionic model. Several ionic models exist and they typically consist of few to hundreds of equations. We remark that there is no restriction in the system \eqref{eq:cbc} for having different ionic models on each cell.

Model \cref{eq:cbc} is a slight simplification of a more detailed model by \citet{TMR21}, where the dynamics at the gap junctions is time dependent and nonlinear in $u_i-u_j$. The simplification adopted here follows from linearization and an equilibrium assumption. This procedure leads to a less computationally intensive model. Solving the complete model and compare the results is subject of a future work.

In this section we adapt the techniques used in \cref{sec:introex} to the full problem \cref{eq:cbc}. 
First, in \cref{sec:spacedisc}, we perform the spatial discretization of the cell-by-cell model, obtaining a differential algebraic equation.
Then, in  \cref{sec:toode}, we reduce the problem to a system of ordinary differential equations on the transmembrane boundary.

We start rewriting \cref{eq:cbc} as follows. 
Let $V_{ij}=\ui-\uj$ be the difference of potential defined on the gap junctions $\Gammaij$ for $1\leq j<i\leq N$. 
Please note that we consider $\Gammaij$ with $j<i$ only, this is to avoid any confusion regarding the sign of $\ui$ and $\uj$ in the definition of $V_{ij}$.
Let $\Gammag = \cup_{1\leq j<i\leq N}\Gammaij$ be the union of all gap junctions, $\Gammaz$ is the transmembrane boundary and $\Gamma=\Gammaz\cup\Gammag$ the union of all internal boundaries. We define $V$ on $\Gamma$ by $V|_{\Gamma_{ij}}=V_{ij}$.
This yields \cref{eq:cbcrw_diff_u_eq_V} instead of \cref{eq:cbc_diff_u_eq_V}.
Condition \eqref{eq:cbc_flux_ue_gammam} yields \cref{eq:cbcrw_flux_uz_gammaz}, while summing \cref{eq:cbc_flux_ue_gammam,eq:cbc_flux_ui_gammaim} we obtain \cref{eq:cbcrw_cont_fluxes} for $j=0$. 
Summing \cref{eq:cbc_permeability} inverting the roles of $i,j$ yields \cref{eq:cbcrw_cont_fluxes} for $j\geq 1$, while taking the difference gives \cref{eq:cbcrw_permeability}.
\begin{subequations}\label{eq:cbcrw}
	\begin{empheq}[left={\empheqbiglbrace~}]{alignat=2}
		-\sigmai \Delta \ui &= 0,
		&& \text{in $\Omegai$ for $i=0,\ldots,N$}, \label{eq:cbcrw_poisson} \\		
		\sigmai \partial_{\bni} \ui +\sigmaj \partial_{\bnj} \uj &=  0,
		&& \text{on $\Gammaij$ for $0\leq j<i\leq N$}, \label{eq:cbcrw_cont_fluxes} \\
		\ui - \uj &=V, && \text{on $\Gammaij$ for $0\leq j< i\leq N$}, \label{eq:cbcrw_diff_u_eq_V}\\
		\sigmaj \partial_{\bnj} \uj - \sigmai \partial_{\bni} \ui &=2 \kappa V,
		&& \text{on $\Gammaij$ for $1\leq j<i\leq N$}, \label{eq:cbcrw_permeability}\\
		\sigmaz \partial_{\bnz} \uz &= \It(V,z),\qquad
		&& \text{on $\Gammaz$},  \label{eq:cbcrw_flux_uz_gammaz}\\				
		\partial_t z &= g(V,z),
		&& \text{on $\Gammaz$}. \label{eq:cbcrw_ionic_var_eq}\\
		\sigmaz \partial_{\bnz} \uz &= 0,\qquad
		&& \text{on $\Sigma$}.  \label{eq:cbcrw_flux_uz_sigma}
	\end{empheq}
\end{subequations}
Model \cref{eq:cbcrw} is equivalent to \cref{eq:cbc}, however it is written in a more ``symmetric'' manner.

\subsection{Spatial discretization of the cell-by-cell model}\label{sec:spacedisc}
We discretize all boundary segments $\Gammaij$ with $M_{ij}$ collocation points $\bx_{ij}^k\in \mathring{\Gamma}_{ij}$, for $k=1,\ldots,M_{ij}$ and $0\leq j<i\leq N$. 
Let $M_i$ be the number of discretization points lying on boundary $\Gammai$, $i=0,\ldots,N$. The total number of collocation points on $\Gamma=\cup_{i=0}^N \Gammai$ is $M=\sum_{0\leq j<i\leq N}M_{ij} = \tfrac{1}{2}\sum_{i=0}^N M_i$. We denote $\bx^l$, $l=1,\ldots,M$, the global collocations points on $\Gamma$ and by $\bx_{i}^k$, $k=1,\ldots,M_i$, the local collocation points on $\Gammai$. 
Note that every $\bx^l$ lies on some $\Gammaij$, hence there are $\bx_{ij}^{k_1}$, $\bx_i^{k_2}$, $\bx_j^{k_3}$ satisfying $\bx^l=\bx_{ij}^{k_1}=\bx_i^{k_2}=\bx_j^{k_3}$.

Let $A_i\in\Rb^{M_i\times M}$ be the boolean connectivity matrix mapping a vector $\bm{v}\in\Rb^M$ of global nodal values on $\Gamma$ to the vector $\bm{v}_i\in\Rb^{M_i}$ of local nodal values on $\Gammai$. Every line of $A_i$ has exactly one non zero element: $(A_i)_{kl}=1$ for $k,l$ such that $\bx^l=\bx_i^k$. Note that $A_i^\top$ maps local to global degrees of freedom. 
We also define $B_i\in\Rb^{M_i\times M}$ having the same sparsity pattern as $A_i$. Let $(B_i)_{kl}$ be the only non zero element in the $k$-th line, hence $\bx^l=\bx_i^k$. If $\bx^l\in\Gammaij$ with $j<i$ then $(B_i)_{kl}=1$, else $(B_i)_{kl}=-1$.

Let $\PSi\in\Rb^{M_i\times M_i}$ be the discrete Poincaré--Steklov operator on each domain $\Omegai$ and $\bui\in\Rb^{M_i}$ the vector of coordinates representing $\gammai_0\ui$. The vector of coordinates $\bV\in\Rb^M$ represents $V$ and $\bVm=A_0\bV$ represents $V|_{\Gamma_0}$. The spatial discretization of \cref{eq:cbcrw_diff_u_eq_V,eq:cbcrw_cont_fluxes} is given by 
\begin{equation}\label{eq:cbcsums}
	\sum_{i=0}^N \sigmai A_i^\top \PSi\bui = \bm{0}, \qquad \qquad
	\sum_{i=0}^N B_i^\top \bui = \bV.
\end{equation}
Recall that $A_0$ is the connectivity matrix mapping a global vector $\bm{v}\in\Rb^M$ to a local vector $\bm{v}_0\in\Rb^{M_0}$ on the transmembrane boundary $\Gamma_0$. Let $M_g=M-M_0$ be the number of points on the gap junctions $\Gamma_g$ and $A_g\in\Rb^{M_g\times M}$ the matrix mapping a global vector to a local vector $\bm{v}_g\in\Rb^{M_g}$ on $\Gamma_g$. 
The spatial discretization of \cref{eq:cbcrw_flux_uz_gammaz,eq:cbcrw_permeability} is
\begin{equation}\label{eq:prepredae}
 \sigma_0\PSi[0]\bm{u}_0 = \It(A_0\bV,\bz),\qquad \qquad 
  \sum_{i=1}^N\sigmai A_g B_i^\top \PSi\bui = -2\kappa A_g\bV.
\end{equation}
As in \cref{sec:introex}, conditions \cref{eq:cbcrw_poisson,eq:cbcrw_flux_uz_sigma} are automatically satisfied by the Green representation formula \eqref{eq:bir} and the definition of the Poincaré--Steklov operator $\PSi[0]$ on $\Omega_0$. Finally, the spatial discretization of \cref{eq:cbcrw_ionic_var_eq} is
\begin{equation}\label{eq:sdgating}
	\bz' = g(A_0\bV,\bz).
\end{equation}
Hence, the spatial discretization of \cref{eq:cbcrw} is given by \cref{eq:cbcsums,eq:prepredae,eq:sdgating}.

\subsection{Reduction to an ordinary differential equation}\label{sec:toode}
In this section we transform the space discretization \cref{eq:cbcsums,eq:prepredae,eq:sdgating} into an ordinary differential equation. First, similarly to \cref{sec:ex_lag}, we search for linear maps
\begin{equation}\label{eq:defphii}
	\mapV_i : \Rb^M\rightarrow \Rb^{M_i}, \qquad\qquad 
	\mapV_i(\bV) = \sigma_i\PSi\bui,\qquad\qquad i=0,\ldots,N,
\end{equation}
where the $\bui$ satisfy \cref{eq:cbcsums}. With the help of these maps we can dispose of  \cref{eq:cbcsums} by inserting \cref{eq:defphii} into \cref{eq:prepredae} and obtain the system of equations
\begin{equation}\label{eq:predae}
	\mapV_0(\bV) = \It(A_0\bV,\bz)=\Cm A_0\bV' + \Iion(A_0\bV,\bz),\qquad \qquad 
	\sum_{i=1}^N A_g B_i^\top \mapV_i(\bV) = -2\kappa A_g\bV.
\end{equation}
However, \cref{eq:predae} is a differential algebraic equation (DAE), which requires more involved time marching schemes than a simple ODE.
Therefore, departing from the definition of the maps $\mapV_i$ given in \cref{thm:phii}, in \cref{thm:ode} we derive a new map which takes into account also the algebraic condition (second equality in \cref{eq:predae}). This new map will allow us to derive an ODE instead of a DAE.

We start with the theorem below, where we compute the maps $\mapV_i$ of \cref{eq:defphii}. The procedure adopted here is inspired from \cite{LaS03}, where a domain decomposition technique for the BEM is presented. 
\begin{theorem}\label{thm:phii}
The linear maps $\mapV_i$ from \cref{eq:defphii} satisfy 
\begin{equation}
	\mapV_i(\bV) = -B_i\blambda,
\end{equation}
with $\blambda\in\Rb^M$ and $\bbeta\in\Rb^N$ solution to
\begin{equation}\label{eq:sys}
	\begin{pmatrix}
		F & G \\ G^\top & 0 
	\end{pmatrix}
	\begin{pmatrix}
		\blambda \\ \bbeta
	\end{pmatrix}
	=
	\begin{pmatrix}
		\bV \\ \bm{0}
	\end{pmatrix}.
\end{equation}
The matrices $F\in\Rb^{M\times M}$, $G\in\Rb^{M\times N}$ are defined by
\begin{equation}
	F = -\sum_{i=0}^N \sigmai^{-1}B_i^\top(\PSi^+)^{-1} B_i, \qquad G=(B_1^\top\be_1,\ldots,B_N^\top\be_N),
\end{equation}
$\be_i\in\Rb^{M_i}$ is the vector of ones and
\begin{equation}
	\PSi^+ = \PSi+\alpha_i \bei\bei^\top,
\end{equation}
with $\alpha_i>0$, $i=0,\ldots,N$. If needed, $\bui$ for $i=0,\ldots,N$ is computed with
\begin{equation}\label{eq:phitoui}
	\bui = -\sigmai^{-1}(\PSi^+)^{-1} B_i\blambda + \beta_i\bei,
\end{equation}
where $\bbeta=(\beta_1,\ldots,\beta_N)^\top$ and $\beta_0=0$.
\end{theorem}
\begin{proof}
	As in \cref{sec:ex_lag} we notice that instead of solving \cref{eq:cbcsums} we can solve a constrained minimization problem with Lagrangian function
	\begin{equation}
		\mathcal{L}(\bui[0],\ldots,\bui[N],\blambda) = \sum_{i=0}^N \frac{\sigmai}{2} \langle \PSi\bui,\bui\rangle + \sum_{i=0}^N\langle B_i^\top\bui,\blambda\rangle-\langle\bV,\blambda\rangle.
	\end{equation}
	Indeed, 
	\begin{equation}\label{eq:nablalambda}
		\nabla_{\blambda}\mathcal{L}(\bui[0],\ldots,\bui[N],\blambda) = \sum_{i=0}^N B_i^\top\bui-\bV = \bm{0}
	\end{equation}
	is equivalent to the second equality in \cref{eq:cbcsums}. The first equality in \cref{eq:cbcsums} follows from
	\begin{equation}\label{eq:nablaui}
		\nabla_{\bui}\mathcal{L}(\bui[0],\ldots,\bui[N],\blambda)=\sigmai\PSi\bui+B_i\blambda = \bm{0}
	\end{equation}
	and 
	\begin{equation}
		\bm{0}=\sum_{i=0}^N A_i^\top \nabla_{\bui}\mathcal{L}(\bui[0],\ldots,\bui[N],\blambda) = \sum_{i=0}^N\sigmai A_i^\top \PSi\bui + A_i^\top B_i\blambda = \sum_{i=0}^N\sigmai A_i^\top \PSi\bui,
	\end{equation}
	where we used $\sum_{i=0}^N A_i^\top B_i=0$ (for every $1$ there is a $-1$). 
	
	Let us solve \cref{eq:nablaui}. We denote $\bei\in\Rb^{M_i}$ the vector of ones, imposing 
	\begin{equation}\label{eq:condlambdaa}
		\langle B_i\blambda,\bei\rangle=0, \qquad i=0,\ldots,N,
	\end{equation}
	a solution $\bui$ to \cref{eq:nablaui} exists and is given by
	\begin{equation}\label{eq:solui}
		\bui = -\sigmai^{-1}(\PSi^+)^{-1} B_i\blambda + \beta_i\bei, \quad \text{with} \quad \PSi^+ = \PSi+\alpha_i \bei\bei^\top,
	\end{equation}
	$\alpha_i>0$ and $\beta_i\in\Rb$. Since if $\bui$, $i=0,\ldots,N$, are solutions to \cref{eq:cbcsums} then also $\bui+C\bei$ are solutions, we choose to set $\beta_0=0$ and remove this degree of freedom.
	Note as well that
	\begin{equation}
		\langle B_0\blambda,\be_0\rangle = -\sum_{i=1}^N\langle B_i\blambda,\bei\rangle \qquad \forall \blambda\in\Rb^{M},
	\end{equation}
	hence \cref{eq:condlambdaa} is replaced by the sufficient one
	\begin{equation}\label{eq:condlambdab}
		\langle B_i\blambda,\bei\rangle=0, \qquad i=1,\ldots,N.
	\end{equation}
	Inserting \cref{eq:solui} into \cref{eq:nablalambda} yields
	\begin{equation}\label{eq:condV}
		\bV=\sum_{i=0}^N B_i^\top(-\sigmai^{-1}(\PSi^+)^{-1} B_i\blambda + \beta_i\bei) = -\sum_{i=0}^N \sigmai^{-1}B_i^\top(\PSi^+)^{-1} B_i\blambda +\sum_{i=1}^N \beta_i B_i^\top\bei,
	\end{equation}
	recall that $\beta_0=0$. From \cref{eq:condV,eq:condlambdab} follow \cref{eq:sys}.
\end{proof}

Now we use the result of \cref{thm:phii} and the second equality of \cref{eq:predae} in order to derive a standard ODE problem. We recall that $\bVm=A_0\bV$.
\begin{theorem}\label{thm:ode}
	The space discretization \cref{eq:cbcsums,eq:prepredae,eq:sdgating}  of \cref{eq:cbc} is equivalent to  the ordinary differential equations system
	\begin{subequations}\label{eq:cbcode}
		\begin{empheq}[left={\empheqbiglbrace~}]{alignat=1}
			\Cm\bVm' + \Iion(\bVm,\bz)&=\mapV(\bVm) ,\\
			\bz' &= g(\bVm,\bz),
		\end{empheq}
	\end{subequations}
where $\mapV(\bVm)=\blambda_0$ and $\blambda_0\in\Rb^{M_0}$, $\blambda_g\in\Rb^{M_g}$, $\bbeta\in\Rb^N$ are solutions to
\begin{equation}\label{eq:odesys}
	\begin{pmatrix}
		F_{00} &  F_{0g} & A_0 G \\
		F_{g0} &  F_{gg}-\kappa^{-1}I & A_g G \\
		G^\top A_0^\top & G^\top A_g^\top & 0
	\end{pmatrix}
\begin{pmatrix}
	\blambda_0 \\ \blambda_g \\ \bbeta 
\end{pmatrix}
=
\begin{pmatrix}
	\bVm \\ \bm{0} \\ \bm{0}
\end{pmatrix},
\end{equation}
with 
\begin{equation}
	F_{00} = A_0FA^\top_0, \quad 
	F_{0g} = A_0FA^\top_g, \quad 
	F_{g0} = A_gFA^\top_0, \quad
	F_{gg} = A_gFA^\top_g.
\end{equation}
\end{theorem}
\begin{proof}
	We denote $\bVg=A_g\bV$, $\bVm=A_0\bV$, $\blambda_g=A_g\blambda$ and $\blambda_0=A_0\blambda$. Also, note that $B_i^\top B_i = A_i^\top A_i$ since $B_i^\top B_i $ projects a global vector forth and back from $\Gammai$ and if a sign change happens it occurs twice. Therefore
	\begin{equation}
		\left(\sum_{i=0}^N B_i^\top B_i\right)\bm{v} = \left(\sum_{i=0}^N A_i^\top A_i \right)\bm{v} = 2 \bm{v} \qquad\forall\, \bm{v}\in\Rb^M,
	\end{equation}
	indeed every segment $\Gammaij$ will receive the contribution from exactly two neighbouring domains.
	From \cref{thm:phii} we have that $\mapV_i(\bV)=-B_i \blambda$, which inserted into the second equality of \cref{eq:predae} yields
	\begin{equation}
	\kappa\bVg = \kappa A_g\bV=-\frac{1}{2}\sum_{i=1}^N A_g B_i^\top \mapV_i(\bV) = \frac{1}{2}A_g\left( \sum_{i=0}^N B_i^\top B_i \right)\blambda = A_g\blambda = \blambda_g.
	\end{equation}
Note as well that $A_0^\top A_0+A_g^\top A_g$ is the identity matrix in $\Rb^M$, hence multiplying the first line $F\blambda+G\bbeta=\bV$ of \cref{eq:sys} with $A_0$ yields
\begin{equation}\label{eq:odesysa}
	\begin{aligned}
	\bVm &=A_0\bV=A_0 F(A_0^\top A_0\blambda+A_g^\top A_g\blambda) + A_0 G\bbeta 
	= F_{00}\blambda_0 + F_{0g}\blambda_g + A_0 G\bbeta.
\end{aligned}
\end{equation}
Similarly, multiplication by $A_g$ yields $\bVg = F_{g0}\blambda_0 + F_{gg}\blambda_g + A_g G\bbeta$ and thus
\begin{equation}\label{eq:odesysb}
	\bm{0} = F_{g0}\blambda_0 + (F_{gg}-\kappa^{-1}I)\blambda_g + A_g G\bbeta .
\end{equation}
For the second line of \cref{eq:sys} we have
\begin{equation}\label{eq:odesysc}
	\bm{0}= G^\top\blambda 
	= G^\top A_0^\top \blambda_0 +  G^\top A_g^\top \blambda_g.
\end{equation}
Relations \cref{eq:odesysa,eq:odesysb,eq:odesysc} yield \cref{eq:odesys}. The identity $-B_0=A_0$ and $\mapV_0(\bV)=-B_0\blambda=A_0\blambda=\blambda_0$ implies $\mapV_0(\bV)=\mapV(\bV_0)$ and hence \cref{eq:cbcode}.
\end{proof}

From the proof of \cref{thm:ode} we see that $\bV$, $\blambda$ of \cref{thm:phii} are given by $\bV = A_0^\top \bV_0+A_g^\top \bVg$, $\blambda = A_0^\top \blambda_0+A_g^\top \blambda_g$, $\kappa\bVg=\blambda_g$ and moreover $\bbeta$ is the same as in \cref{thm:phii}; hence, if needed, $\bui$ for $i=0,\ldots,N$ can be computed as in \cref{thm:phii}. Note as well that in \cref{eq:odesys} we have chosen to use $\blambda_g$ as unknown, instead of the alternative $\bVg$. If we used $\bVg$ we would obtain the same matrix as in \cref{eq:odesys} but with the second column multiplied by $\kappa$ and therefore break the symmetry. 

\subsection{Time integration}

For the time integration of \cref{eq:cbcode}, we use the multirate explicit stabilized method mRKC \cite{AGR22} for problems
\begin{equation}\label{eq:mode}
	\bm{y}'=f_F(t,\bm{y})+f_S(t,\bm{y}),\qquad \bm{y}(0)=\bm{y}_0,
\end{equation}
where $f_F$ is a stiff term and $f_S$ is a mildly stiff but more expensive term. The mRKC scheme is fully explicit and does not have any step size restriction. Its stability properties are inherited from the RKC methods \cite{HoS80}, which use an increased number of stages, with respect to classical methods, to increase stability. Since stability grows quadratically with the work load, the methods are particularly efficient.
For the integration of \cref{eq:cbcode} with mRKC, we rewrite \cref{eq:cbcode} as \cref{eq:mode}, with
\begin{equation}
	\bm{y} = \begin{pmatrix}
		\bVm \\ \bz
	\end{pmatrix},\qquad
	f_F(t,\bm{y})=\begin{pmatrix}
		\mapV(\bVm)/\Cm \\ \bm{0}
	\end{pmatrix},\qquad 
	f_S(t,\bm{y})=\begin{pmatrix}
		-(\Iion(\bVm,\bz)+\Istim(t))/\Cm \\ g(\bVm,\bz)
	\end{pmatrix},
\end{equation} 
where $\Istim(t)$ is used to stimulate some cells and initiate an action potential propagation.

\section{Numerical experiments}\label{sec:num}
In this section we perform some numerical experiments in order to asses the accuracy of the space-time discretization of the cell-by-cell model \eqref{eq:cbc} but also investigate the regularity properties of the model itself.

We start with two experiments, in \cref{sec:conv_psi,sec:imp_dic_par_cv}, where we investigate the convergence rates of the maps $\mapV_i$ from \cref{thm:phii} and then the impact of the mesh and step size on the accuracy of the conduction velocity (CV). These experiments are crucial to understand which discretization parameters yield solutions within a certain error tolerance. 

In the subsequent experiments the goal is to study the model itself. For instance, in \cref{sec:dep_cv_kappa} we investigate the effect of the gap junction's permeability $\kappa$ on the CV and in \cref{sec:dep_cv_area} we study how the contact area between cells affects CV.
Before presenting the results, we resume here below our computational setting.

\paragraph{Computational setup} The following numerical experiments have been performed with our C++ code, where for the dense linear algebra routines we employ the \texttt{Eigen} library \cite{GuB10}. The ionic model, is taken from~\texttt{CellML} \cite{CLN03} and the relative C code is produced with the \texttt{Myokit} library \cite{CCL16}.  Concerning the model \cref{eq:cbc}, the number of cells $N$ and the domains $\Omega_i$ vary from one experiment to another and are specified later. If not stated otherwise, in the next experiments we use the coefficients $\Cm$, $\sigmai$, $\kappa$, given in \cref{tab:model_coeff}. The values for $\Cm$, $\sigmai$, $\kappa$ are taken from \cite{SHL05}, where for $\kappa$ we consider $\kappa=1/R_m$ with $R_m=0.00145\ \si{\kilo\ohm\square\centi\metre}$. If not specified, we consider the ionic model from \citet{CRN98}. The initial values for $V$ for and the ionic model's state variables are uniform on the transmembrane boundary and are taken from the \texttt{Myokit}'s code. For instance, for the Courtemanche-Nattel-Ramir\'ez model the initial value for $V$ is $V_0=-81.18\ \mV$. 

\begin{table}[htb]
	\begin{center}
		\begin{tabular}{c|c|c|c}
			$C_m$ & $\sigmae$ & $\sigma_1,\ldots,\sigma_N$ & $\kappa$ \\ \hline
			\SI{1}{\micro\farad\per\square\centi\metre} & \SI{20}{\milli\siemens\per\centi\metre} & \SI{3}{\milli\siemens\per\centi\metre} & \SI{690}{\milli\siemens\per\square\centi\metre}
		\end{tabular}
	\end{center}
	\caption{Model's coefficients employed in numerical experiments.}
	\label{tab:model_coeff}
\end{table}

\subsection{Approximation properties of the $\mapV$ operators}\label{sec:conv_psi}
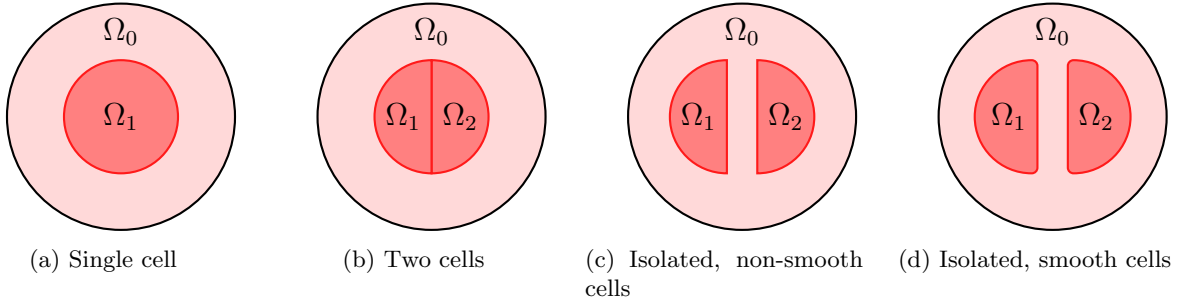
\begin{figure}[t]
	\begin{center}
		\begin{subfigure}[t]{0.23\textwidth}
			\centering
			\includestandalone{../images/conv_psi/illustrations/circle}
			\caption{Single cell}
			\label{fig:ill_conv_psi_a}
		\end{subfigure}
		\hfill
		\begin{subfigure}[t]{0.23\textwidth}
			\centering
				
\pgfmathsetmacro{\cd}{0.0}
\pgfmathsetmacro{\crz}{1.5}
\pgfmathsetmacro{\cru}{0.75}
\pgfmathsetmacro{\crd}{0.0}
\begin{tikzpicture}	
	\fill[cellexterior,draw=black,thick] (0,0) circle (\crz);		
	
	\filldraw[fill=cellinterior,draw=cellmembrane,thick]
	(\cd,0) -- (\cd,{-\cru+\crd}) arc (180:270:\crd) arc(-90:90:\cru) arc(90:180:\crd) -- (\cd,0);
	
	\filldraw[fill=cellinterior,draw=cellmembrane,thick]
	(-\cd,0) -- (-\cd,{-\cru+\crd}) arc (0:-90:\crd) arc(-90:-270:\cru) arc(90:0:\crd) -- (-\cd,0);
	
	\node[font=\large] at ({-\cd-\cru/2},0) {$\Omega_1$};
	\node[font=\large] at ({\cd+\cru/2},0) {$\Omega_2$};
	\node[font=\large] at (0,{0.5*(\cru+\crz)}) {$\Omega_0$};

\end{tikzpicture}
	
			\caption{Two cells}
			\label{fig:ill_conv_psi_b}
		\end{subfigure}
		\hfill
		\begin{subfigure}[t]{0.23\textwidth}
			\centering
				
\pgfmathsetmacro{\cd}{0.2}
\pgfmathsetmacro{\crz}{1.5}
\pgfmathsetmacro{\cru}{0.75}
\pgfmathsetmacro{\crd}{0.0}
\begin{tikzpicture}	
	\fill[cellexterior,draw=black,thick] (0,0) circle (\crz);		
	
	\filldraw[fill=cellinterior,draw=cellmembrane,thick]
	(\cd,0) -- (\cd,{-\cru+\crd}) arc (180:270:\crd) arc(-90:90:\cru) arc(90:180:\crd) -- (\cd,0);
	
	\filldraw[fill=cellinterior,draw=cellmembrane,thick]
	(-\cd,0) -- (-\cd,{-\cru+\crd}) arc (0:-90:\crd) arc(-90:-270:\cru) arc(90:0:\crd) -- (-\cd,0);
	
	\node[font=\large] at ({-\cd-\cru/2},0) {$\Omega_1$};
	\node[font=\large] at ({\cd+\cru/2},0) {$\Omega_2$};
	\node[font=\large] at (0,{0.5*(\cru+\crz)}) {$\Omega_0$};

\end{tikzpicture}
	
			\caption{Isolated, non-smooth cells}
			\label{fig:ill_conv_psi_c}
		\end{subfigure}
		\hfill
		\begin{subfigure}[t]{0.23\textwidth}
			\centering
				
\pgfmathsetmacro{\cd}{0.2}
\pgfmathsetmacro{\crz}{1.5}
\pgfmathsetmacro{\cru}{0.75}
\pgfmathsetmacro{\crd}{0.075}
\begin{tikzpicture}	
	\fill[cellexterior,draw=black,thick] (0,0) circle (\crz);		
	
	\filldraw[fill=cellinterior,draw=cellmembrane,thick]
	(\cd,0) -- (\cd,{-\cru+\crd}) arc (180:270:\crd) arc(-90:90:\cru) arc(90:180:\crd) -- (\cd,0);
	
	\filldraw[fill=cellinterior,draw=cellmembrane,thick]
	(-\cd,0) -- (-\cd,{-\cru+\crd}) arc (0:-90:\crd) arc(-90:-270:\cru) arc(90:0:\crd) -- (-\cd,0);
	
	\node[font=\large] at ({-\cd-\cru/2},0) {$\Omega_1$};
	\node[font=\large] at ({\cd+\cru/2},0) {$\Omega_2$};
	\node[font=\large] at (0,{0.5*(\cru+\crz)}) {$\Omega_0$};

\end{tikzpicture}
	
			\caption{Isolated, smooth cells}
			\label{fig:ill_conv_psi_d}
		\end{subfigure}	
	\end{center}
	\caption{Illustration of the geometrical settings employed in \cref{sec:conv_psi}.}
	\label{fig:ill_conv_psi}
\end{figure}

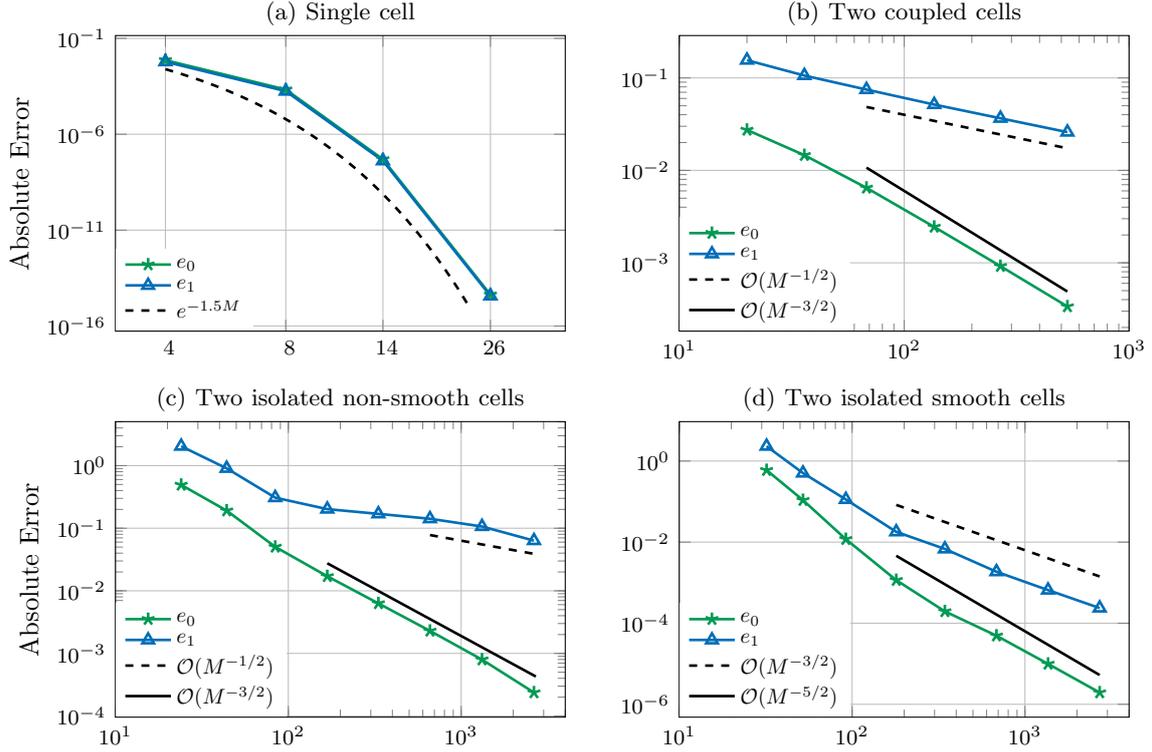
\begin{figure}[t]
\centering
\begin{tikzpicture}

\pgfplotsset{
    log x ticks with fixed point/.style={
    xticklabel={
      \pgfkeys{/pgf/fpu=true}
      \pgfmathparse{exp(\tick)}%
      \pgfmathprintnumber[fixed relative, precision=3]{\pgfmathresult}
      \pgfkeys{/pgf/fpu=false}
    }
  },
}

\begin{groupplot}[
    group style = {
        group size=2 by 2,
        horizontal sep=1.5cm,
        vertical sep=1.2cm,
        ylabels at = edge left,
    },
    grid = major,
    width = 7.5cm,
    height = 5.5cm,
    title style = {text width = 6cm},
    legend cell align=left,
    legend pos = south west,
    legend style={nodes={scale=0.8, transform shape}},
    legend style={draw=\legendboxdraw,fill=white,at={(0,0)}},
    ymode=log,xmode=log,
    title style={anchor=north, yshift=4ex},
    every axis plot/.append style=
    {
        solid,
        line width=\plotlinewidth pt,
        mark size=\plotmarksizeu pt,
    },
    tick label style = {font=\footnotesize},
    ylabel = {Absolute Error},
    xmin=10,xmax=4e3,
]

\nextgroupplot[title={\subcaption{Single cell\label{fig:plot_conv_psi_a}}},xmin=3,xmax=40,log x ticks with fixed point,xtick = data]

\addplot[colorone,mark=\markone] table [x=n_dofs,y=err_trace,col sep=comma] {\currfiledir conv_test_1_mdom_0.csv};
\addlegendentry{$e_0$}
\addplot[colortwo,mark=\marktwo] table [x=n_dofs,y=err_deriv,col sep=comma] {\currfiledir conv_test_1_mdom_0.csv};
\addlegendentry{$e_1$}
\addplot[black,dashed,domain=4:23] {exp(-1.5*x)};
\addlegendentry{$e^{-1.5 M}$}

\nextgroupplot[title={\subcaption{Two coupled cells\label{fig:plot_conv_psi_b}}},xmax=1e3]

\addplot[color=colorone,mark=\markone] table [x=n_dofs,y=err_trace,col sep=comma] {\currfiledir conv_test_1_mdom_-1.csv};
\addlegendentry{$e_0$}

\addplot[color=colortwo,mark=\marktwo] table [x=n_dofs,y=err_deriv,col sep=comma] {\currfiledir conv_test_1_mdom_-1.csv};
\addlegendentry{$e_1$}

\addplot[black,dashed,domain=68:532] {0.4*(x^(-0.5))};
\addlegendentry{$\mathcal{O}(M^{-1/2})$}
\addplot[black,solid,domain=68:532] {6*(x^(-1.5))};
\addlegendentry{$\mathcal{O}(M^{-3/2})$}

\nextgroupplot[title={\subcaption{Two isolated non-smooth cells\label{fig:plot_conv_psi_c}}}]

\addplot[color=colorone,solid,line width=\plotlinewidth pt,mark=\markone,mark size=\plotmarksizeu pt] table [x=n_dofs,y=err_trace,col sep=comma] 
{\currfiledir conv_test_2_mdom_-2.csv}; \addlegendentry{$e_0$}
\addplot[color=colortwo,line width=\plotlinewidth pt,mark=\marktwo,mark size=\plotmarksized pt] table [x=n_dofs,y=err_deriv,col sep=comma] 
{\currfiledir conv_test_2_mdom_-2.csv}; \addlegendentry{$e_1$}
\addplot[black,dashed,domain=660:2704] {2*(x^(-0.5))} ; \addlegendentry{$\mathcal{O}(M^{-1/2})$}
\addplot[black,solid,domain=168:2704] {6e1*(x^(-1.5))} ; \addlegendentry{$\mathcal{O}(M^{-3/2})$}

\nextgroupplot[title={\subcaption{Two isolated smooth cells\label{fig:plot_conv_psi_d}}}]

\addplot[color=colorone,solid,line width=\plotlinewidth pt,mark=\markone,mark size=\plotmarksizeu pt] table [x=n_dofs,y=err_trace,col sep=comma] 
{\currfiledir conv_test_2_mdom_-5.csv}; \addlegendentry{$e_0$}
\addplot[color=colortwo,line width=\plotlinewidth pt,mark=\marktwo,mark size=\plotmarksized pt] table [x=n_dofs,y=err_deriv,col sep=comma] 
{\currfiledir conv_test_2_mdom_-5.csv}; \addlegendentry{$e_1$}
\addplot[black,dashed,domain=180:2704] {2e2*(x^(-1.5))} ; \addlegendentry{$\mathcal{O}(M^{-3/2})$}
\addplot[black,solid,domain=180:2704] {2e3*(x^(-2.5))} ; \addlegendentry{$\mathcal{O}(M^{-5/2})$}

\end{groupplot}

\end{tikzpicture}
	
\caption{Convergence rates of the $\mapV$ operator defined in \cref{thm:phii} for the problems depicted in \cref{fig:ill_conv_psi}.}
\label{fig:plot_conv_psi}
\end{figure}

In this experiment we investigate the convergence rate of the $\mapV_i$ operators defined in \cref{thm:phii}. For this purpose, we conduct four convergence experiments, one for every geometry displayed in \cref{fig:ill_conv_psi}, and we display the errors on the trace and normal derivative against the number of degrees of freedom $M$ in \cref{fig:plot_conv_psi}.

Let us describe the geometries of \cref{fig:ill_conv_psi}. In \cref{fig:ill_conv_psi_a} we have a model composed of one cell ($N=1$), defined by $\Omega_0 = \{\bx\in\Rb^2 : 2< \Vert\bx\Vert < 4\}$ and $\Omega_1 = \{\bx\in\Rb^2 : \Vert\bx\Vert < 2\}$. In \cref{fig:ill_conv_psi_b} we split the cell by introducing a vertical gap junction, hence we have the same $\Omega_0$ but $\Omega_1 = \{\bx\in\Rb^2 : \Vert\bx\Vert<2, x_1<0\}$ and $\Omega_2 = \{\bx\in\Rb^2 : \Vert\bx\Vert<2, x_1>0\}$. In \cref{fig:ill_conv_psi_c} we keep the same cells but remove the gap junction by introducing an horizontal gap of size $0.4$ between $\Omega_1$ and $\Omega_2$. Finally, in \cref{fig:ill_conv_psi_d} we keep the separation of $\Omega_1$, $\Omega_2$ but smooth out the corners by introducing quarter of circles of radius \SI{0.2}{\um}.

For the setting of \cref{fig:ill_conv_psi_a,fig:ill_conv_psi_b} an exact solution to \cref{eq:cbcrw_poisson,eq:cbcrw_diff_u_eq_V,eq:cbcrw_cont_fluxes,eq:cbcrw_permeability}, with $V$ defined by \cref{eq:cbcrw_diff_u_eq_V}, is given by
\begin{equation}\label{eq:defu0u1}
	\begin{aligned}
		u_0(\bx) &= \frac{\sigma_1}{\sigma_0} \frac{16+\Vert \bx\Vert^2}{6\Vert \bx\Vert^2}x_2, \qquad & u_1(\bx) = u_2(\bx)&= -\frac{1}{2}x_2.
	\end{aligned}
\end{equation}
Therefore, for different values of $M$ (i.e.\ number of collocation points), we can compute the vector of coefficients $\bui$, define $\bV$ as in \cref{eq:cbcsums}, solve \cref{eq:sys} and compute the errors
\begin{equation}\label{eq:def_err01}
\begin{aligned}
	e_1 &=\max_{i=0,\ldots,N} \Vert \mapV_i(\bV)-\sigma_i\partial_{\bn_i}u_i\Vert_{L^2(\Gamma_i)}, \qquad & 
	e_0 &= \max_{i=0,\ldots,N} \Vert \widetilde\mapV_i(\bV)-u_i\Vert_{L^2(\Gamma_i)/\Rb},
\end{aligned}
\end{equation}
where $\widetilde\mapV_i(\bV)=-\sigmai^{-1}(\PSi^+)^{-1} B_i\blambda + \beta_i\bei$ and thus $\widetilde\mapV_i(\bV_m)$ approximates $u_i$ (cf. \cref{eq:phitoui}), up to a constant.

For the geometries of \cref{fig:ill_conv_psi_c,fig:ill_conv_psi_d} we do not possess an exact solution. 
Hence, we set $V(\bx)=\cos(\pi x_1)\sin(\pi x_2)$ and errors $e_1$, $e_0$ are now computed as
\begin{equation}\label{eq:def_err01_ref}
	\begin{aligned}
		e_1 &=\max_{i=0,\ldots,N} \Vert \mapV_i(\bV)-\mapV_i(\bV^*)\Vert_{L^2(\Gamma_i)}, \qquad & 
		e_0 &= \max_{i=0,\ldots,N} \Vert \widetilde\mapV_i(\bV)- \widetilde\mapV_i(\bV^*)\Vert_{L^2(\Gamma_i)/\Rb},
	\end{aligned}
\end{equation}
where $\mapV_i(\bV^*)$, $\widetilde\mapV_i(\bV^*)$ are reference solutions calculated on a finer mesh.

We display the errors $e_0,e_1$ with respect to $M$, for the geometries of \cref{fig:ill_conv_psi}, in \cref{fig:plot_conv_psi}. Due to the smoothness of the solutions and the boundaries, we remark in \cref{fig:plot_conv_psi_a} that for the first problem we obtain exponential convergence thanks to the trigonometric Lagrange basis functions; with very few degrees of freedom machine precision is achieved. This result is in line with the theory and experiments performed in \cite{HeJ18}. 
In the second problem, the boundary is Lipschitz continuous only, which prevents exponential convergence. Indeed, in \cref{fig:plot_conv_psi_b} the convergence rates for the trace and normal derivative are $1.5$ and $0.5$, respectively.

In \cref{fig:plot_conv_psi_a,fig:plot_conv_psi_b} we see how the convergence rates decrease dramatically when the circle is divided into two half-circles. The purpose of two last experiments is to demonstrate numerically that this phenomenon is due to the non smooth boundaries, rather than the introduction of a gap junction. Indeed, in \cref{fig:plot_conv_psi_c} we observe the same convergence rates as in \cref{fig:plot_conv_psi_b}, while in \cref{fig:plot_conv_psi_d} we obtain higher convergence rates.

\subsection{Impact of discretization parameters on conduction velocity}\label{sec:imp_dic_par_cv}
It is already known that discretization methods, mesh size and step size affect conduction velocity (CV) in the monodomain and bidomain models for cardiac electrophysiology \cite{CBC11,NKB11,PHS16}. In this experiment we investigate how mesh and step size affect the CV for the cell-by-cell model discretized with the BEM in space and the mRKC method \cite{AGR22} in time. In order to be able to employ relatively uniform mesh sizes in this experiment we consider rectangular cells.

\begin{figure}[t]
	\begin{center}
		\includegraphics[width=\textwidth]
		{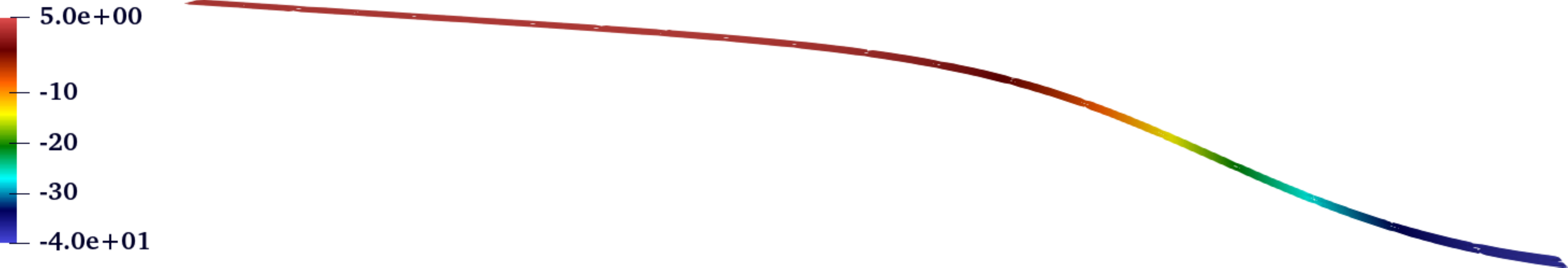}
	\end{center}
	\caption{Action potential propagation in an array of $2\times 20$ cells ($t=2\ms$).}
	\label{fig:arraycells}
\end{figure}

To measure the CV we design the following experiment.
We consider an array of $2\times 30$ connected rectangular cells of width $c_w=\SI{20}{\um}$ and length $c_l=\SI{100}{\um}$, cells are positioned so that their bottom left vertex has coordinates $(i\cdot c_l,j\cdot c_w)$ for $i=0,\ldots,29$, $j=0,1$; yielding a block of cells of width $2\cdot c_w=\SI{40}{\um}$ and length $30\cdot c_l=\SI{3000}{\micro\metre}$. The outer domain is a rectangle of size $\SI{440}{\um}\times \SI{5000}{\um}$ centred on the array of cells. 
To initiate an action potential traversing the cell's array a stimulus of $\SI{300}{\micro\ampere\per\square\centi\metre}$ is applied for a duration of \SI{1}{\milli\second} at the transmembrane boundary of the two leftmost cells. CV is computed as the average over $CV_k=\Vert \bm{p}_k-\bm{q}_k\Vert/(t_{\bm{p}_k}-t_{\bm{q}_k})$ for $k=1,\ldots,5$, where: $\bm{p}_k=((7.5+k)\cdot c_l)\,\um$, $\bm{q}_k=((17.5+k)\cdot c_l,0)\,\um$ and $t_{\bm{p}_k}$, $t_{\bm{q}_k}$ are the time instants in which $V$ exceeds the threshold of $V_\mathrm{th}=\SI{-20}{\mV}$ in $\bm{p}_k$, $\bm{q}_k$, respectively. The choice of $\bm{p}_k$, $\bm{q}_k$ is such that measures are taken sufficiently far from the stimulated point and to avoid boundary effects as well.
See \cref{fig:arraycells} for an illustration of the solution at time $t=\SI{2}{\ms}$.

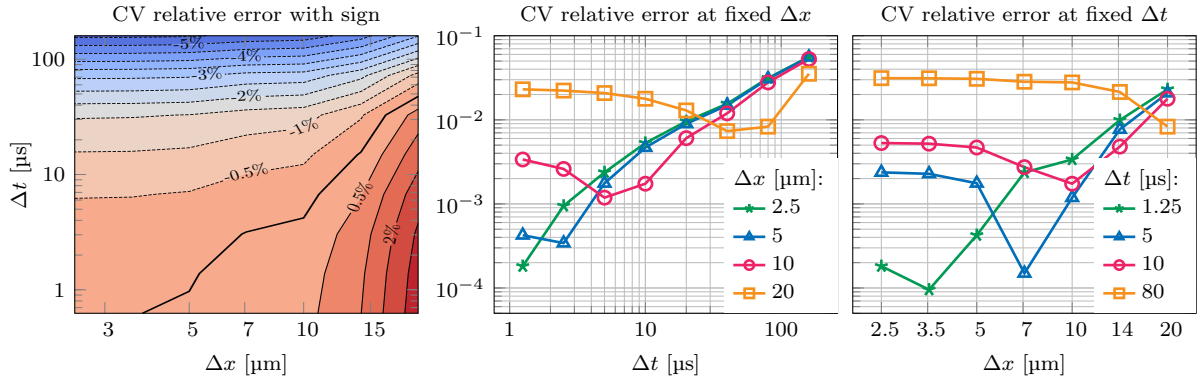
\begin{figure}[t]
	\centering
	\begin{tikzpicture}[font=\footnotesize]

\pgfplotscreateplotcyclelist{mycolorlisttest}{%
	color=colorone,mark=\markone\\
	color=colortwo,mark=\marktwo\\
	color=colorthree,mark=\markthree\\
	color=colorfour,mark=\markfour\\
	color=colorfive,mark=\markfive\\
	color=colorsix,mark=\marksix\\
}

\pgfplotsset{
    log x ticks with fixed point/.style={
    xticklabel={
      \pgfkeys{/pgf/fpu=true}
      \pgfmathparse{exp(\tick)}%
      \pgfmathprintnumber[fixed relative, precision=3]{\pgfmathresult}
      \pgfkeys{/pgf/fpu=false}
    }
  },
}

\begin{groupplot}[
	group style = {
		group size=3 by 1,
	},
    every axis plot/.append style=
    {
        solid,
        line width=\plotlinewidth pt,
        mark size=\plotmarksizeu pt,
    },
	cycle list name = mycolorlisttest,
    width = 6.1cm,
	legend cell align=left,legend pos = south east,
    legend style={
		font=\footnotesize,
		draw=\legendboxdraw,fill=white,at={(1,0)},
		anchor=south east},
	legend image post style={scale=0.7},
	title style={anchor=north, yshift=2ex},
	log x ticks with fixed point
]

\nextgroupplot[
	title = {CV relative error with sign},
	xlabel={$\Delta x\ [\um]$},
	ylabel={$\Delta t\ [\si{\micro\second}]$},
	ylabel shift={-7pt},
	ymode = log, xmode = log,
	ymin = 0.625, ymax = 160,
	xmin = 2.5, xmax = 20,
	xtick = {3,5,7,10,15},
	enlargelimits=false,
	axis on top=true,
	log ticks with fixed point,
]

\addplot graphics[xmin=2.5,xmax=20,ymin=0.625,ymax=160] {\currfiledir CVcontour};


\nextgroupplot[title={CV relative error at fixed $\Delta x$},
xlabel={$\Delta t$ [\si{\micro\second}]},
xtick={0.001,0.01,0.1},
xticklabels={1,10,100},
grid=both,xmode=log,ymode=log,ymax=1e-1,ymin=5e-5]

\addlegendimage{empty legend}
\addlegendentry{\hspace{-.5cm}$\Delta x\ [\um]$:}
\addplot+ table [x=dt,y=absCVerr_dG_1,col sep=comma] {\currfiledir CV_vs_dt.csv};
\addlegendentry{$2.5$}
\addplot+ table [x=dt,y=absCVerr_dG_3,col sep=comma] {\currfiledir CV_vs_dt.csv};
\addlegendentry{$5$}
\addplot+ table [x=dt,y=absCVerr_dG_5,col sep=comma] {\currfiledir CV_vs_dt.csv};
\addlegendentry{$10$}
\addplot+ table [x=dt,y=absCVerr_dG_7,col sep=comma] {\currfiledir CV_vs_dt.csv};
\addlegendentry{$20$}


\nextgroupplot[title={CV relative error at fixed $\Delta t$},
grid=both,xshift=-0.8cm,xlabel={$\Delta x$ [\si{\micro\metre}]},
xtick={2.5,3.5,5,7,10,14,20},
log ticks with fixed point,
yticklabels=\empty,xmode=log,ymode=log,ymax=1e-1,ymin=5e-5]

\addlegendimage{empty legend}
\addlegendentry{\hspace{-0.5cm}$\Delta t\ [\si{\micro\second}]$:}
\addplot+ table [x=dx,y=absCVerr_dt_1,col sep=comma] {\currfiledir CV_vs_dx.csv}; 
\addlegendentry{$1.25$}
\addplot+ table [x=dx,y=absCVerr_dt_3,col sep=comma] {\currfiledir CV_vs_dx.csv}; 
\addlegendentry{$5$}
\addplot+ table [x=dx,y=absCVerr_dt_4,col sep=comma] {\currfiledir CV_vs_dx.csv}; 
\addlegendentry{$10$}
\addplot+ table [x=dx,y=absCVerr_dt_7,col sep=comma] {\currfiledir CV_vs_dx.csv}; 
\addlegendentry{$80$}

\end{groupplot}
\end{tikzpicture}
	\caption{Conduction velocity accuracy with respect to timestep $\Delta t$ and mesh size $\Delta x$.}
	\label{fig:cv_vs_dy_dx_dt}
\end{figure}

We solve \cref{eq:cbcode} with different step sizes $\Delta t$ and mesh size $\Delta x$ and compute the signed relative error on CV:
$E_\mathrm{CV} = (\mathrm{CV}-\mathrm{CV}^*)/\mathrm{CV}^*$,
with $\mathrm{CV}^*$ a reference solution. We display $E_\mathrm{CV}$ as function of $\Delta t$, $\Delta x$ in \cref{fig:cv_vs_dy_dx_dt}. The reference value of CV is $\mathrm{CV}^* \approx 1.27153$.

First, we notice that for coarse space grids the true $\mathrm{CV}^*$ tends to be overestimated, whereas for large time steps it is underestimated. Then, we remark that even with relatively large mesh sizes $\Delta x=\SI{20}{\um}$ the estimated CV remains within a $2\%$ error.
\Cref{fig:cv_vs_dy_dx_dt} (middle and right panels) also shows the same results, but for the unsigned relative error $\vert E_\mathrm{CV} \vert$ and fixing either $\Delta x$ or $\Delta t$. We observe that the local minimal appearing in the curves is due to the cancellation of the positive spatial discretization error with the negative time discretization error.

Based on the results of this section, in the forthcoming experiments we consider $\Delta t\leq \SI{0.02}{\ms}$ and $\Delta x\leq \SI{10}{\um}$, which, for this experiment, yield a relative error of less than $5\%$.

\subsection{Dependence of conduction velocity on gap junctions' permeability and cells inner conductivity}\label{sec:dep_cv_kappa}

\begin{figure}[t]
\centering
\begin{tikzpicture}

\pgfplotscreateplotcyclelist{mycolorlisttest}{%
color=colorone,mark=\markone\\
color=colortwo,mark=\marktwo\\
color=colorthree,mark=\markthree\\
color=colorfour,mark=\markfour\\
color=colorfive,mark=\markfive\\
color=colorsix,mark=\marksix\\
}

\pgfplotsset{
    log x ticks with fixed point/.style={
    xticklabel={
      \pgfkeys{/pgf/fpu=true}
      \pgfmathparse{exp(\tick)}%
      \pgfmathprintnumber[fixed relative, precision=3]{\pgfmathresult}
      \pgfkeys{/pgf/fpu=false}
    }
  },
}

\begin{groupplot}[
        group style = {
            group size=2 by 1,
            horizontal sep=1.5cm,
        },
        grid = major,
        width = 7.5cm,
        height = 5.5cm,
        title style = {text width = 6cm},
        cycle list name = mycolorlisttest,
        legend cell align=left,
        legend pos = north west,
        xmode=log,
        title style={anchor=north, yshift=4ex},
        every axis plot/.append style=
        {
            solid,
            line width=\plotlinewidth pt,
            mark size=\plotmarksizeu pt,
        },
        tick label style = {font=\footnotesize},
        ylabel = {CV [\si{\meter\per\second}]},
    ]
    
    \nextgroupplot[
        title={\subcaption{CV as function of permeability $\kappa$\label{fig:cv_vs_kappa}}},
        xlabel={$\kappa$ [\si{\milli\siemens\per\square\centi\metre}]},
    ]
    
    \addplot+ table [x=k,y=CV,col sep=comma] {\currfiledir CVl_C_Rl_cl_100.csv};
    \addlegendentry{$c_l=\SI{100}{\um}$}
    \addplot+ table [x=k,y=CV,col sep=comma] {\currfiledir CVl_C_Rl_cl_50.csv};
    \addlegendentry{$c_l=\SI{50}{\um}$}

    \nextgroupplot[
        title={\subcaption{CV as function of conductivity $\sigma_i$\label{fig:cv_vs_si}}},
        extra x ticks = {3},
        extra x tick labels = {3},
        extra x tick style = {color=red,major grid style=red,},
    ]
    
    \addplot+ table [x=si,y=CV,col sep=comma] {\currfiledir CVl_C_si_cl_100.csv};
    \addlegendentry{$c_l=\SI{100}{\um}$}
    \addplot+ table [x=si,y=CV,col sep=comma] {\currfiledir CVl_C_si_cl_50.csv};
    \addlegendentry{$c_l=\SI{50}{\um}$}

\end{groupplot}
\end{tikzpicture}
	
\caption{Dependence of CV on permeability $\kappa$ and intracellular conductivity $\sigma_i$ (red: reference value).}
\label{fig:cv_ap_vs_kappa}
\end{figure}
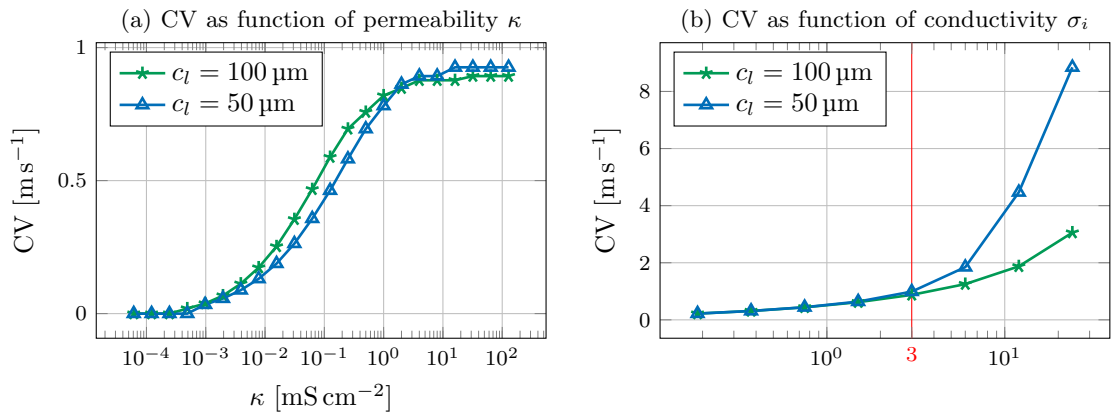

In this experiment we study how CV depends on the gap junctions' permeability $\kappa$ and the inner conductivity $\sigma_i$, $i=1,\ldots,N$.
For that purpose, we consider again an array of $2\times 30$ cells and cells of size $c_w\times c_l$, with fixed $c_w=\SI{10}{\um}$ and either $c_l=\SI{100}{\um}$ or $c_l=\SI{50}{\um}$.
First, we measure CV for varying $\kappa$ but keeping the other coefficients fixed, a stimulus is initiated applying a stimulus of $\SI{200}{\micro\ampere\per\square\centi\metre}$ to the transmembrane boundary of the two leftmost cells. Results for $c_l=\SI{100}{\um}$ and $c_l=\SI{50}{\um}$ are displayed in \cref{fig:cv_vs_kappa}. We observe as CV decreases with $\kappa$ and also that the physiological value $\kappa=690$ is in the range where CV is maximal. We note that for values of $\kappa\leq 2\cdot 10^{-4}$ the action potential does not propagate.
Then, we measure CV for varying $\sigma_i$, $i=1,\ldots,N$, and fixed $\kappa$, results are displayed in \cref{fig:cv_vs_si} for $c_l=\SI{100}{\um}$ and $c_l=\SI{50}{\um}$. Conduction velocity increases with $\sigma_i$, specially for the shorter cells.

\subsection{Dependence of conduction velocity on gap junctions' surface area}\label{sec:dep_cv_area}

\begin{figure}[t]
	\begin{center}
		\begin{subfigure}[t]{\subfigsized\textwidth}
			\centering
				\begin{tikzpicture}
		\begin{scope}[on background layer={color=cellinterior}]
			\fill (-2,0) rectangle (4,2);
		\end{scope}
	
		\draw[draw=cellmembrane, thick] (-2,2)--(4,2);
		\draw[draw=cellmembrane, thick] (-2,0)--(4,0);
		\draw [draw=cellmembrane, thick,  domain=0:1, samples=100] 
		plot ( {-1*(1)*(-1+cos(3*360*\x))}, {2*\x} );
		
		
		\draw[thick,black,<->] (1,1.5)--(2,1.5) ;
		\node [black] at (1.5,1.35) {$a$};
		\node [black] at (3,1.5) {{$k=3$}};
		
		\node [black] at (-1,0.5) {{$\Omega_i$}};
		\node [black] at (3,0.5) {{$\Omega_{i+1}$}};
				
	\end{tikzpicture}
			\caption{Two cells with an intercalated discs at gap junctions.}
			\label{fig:ill_intercalated_discs}
		\end{subfigure}
		\\[1em]
		\begin{subfigure}[t]{\subfigsized\textwidth}
			\centering
				
	\begin{tikzpicture}
		\begin{axis}[height=\aspectratio*\plotsized,width=\plotsized,legend columns=1,legend style={draw=\legendboxdraw,fill=\legendboxfill,at={(0,0)},anchor=south west},legend cell align={left},
			xlabel={$k$}, ylabel={CV [\si{\metre\per\second}]},label style={font=\normalsize},tick label style={font=\normalsize},legend image post style={scale=\legendmarkscale},legend style={nodes={scale=\legendfontscale, transform shape}},
			width = 6.5cm,
			height = 5cm,
			grid = both,
			tick label style = {font=\footnotesize},
			]
			\addplot[color=colorone,solid,line width=\plotlinewidth pt,mark=\markone,mark size=\plotmarksizeu pt] table [x=freq,y=CV,col sep=comma] 
			{\currfiledir CV_vs_freq.csv};
		\end{axis}
	\end{tikzpicture}
	
			\caption{Impact of intercalated discs frequency on CV for $a=\SI{0.5}{\um}$.}
			\label{fig:plot_cv_vs_freq}
		\end{subfigure}\hfill
		\begin{subfigure}[t]{\subfigsized\textwidth}
			\centering
				
	\begin{tikzpicture}
		\begin{axis}[
			xlabel={$a$ [\um]},
			ylabel={CV [\si{\metre\per\second}]},
			label style={font=\normalsize},
			tick label style={font=\normalsize},
			legend image post style={scale=\legendmarkscale},
			xtick={0,1e-4,2e-4},
			xticklabels={0,1,2},
			scaled x ticks = false,
			width = 6.5cm,
			height = 5cm,
			grid = both,
			tick label style = {font=\footnotesize},
			try min ticks=5,
			]
			\addplot[color=colorone,solid,line width=\plotlinewidth pt,mark=\markone,mark size=\plotmarksizeu pt] table [x=trueamp,y=CV,col sep=comma] 
			{\currfiledir CV_vs_amp.csv};
		\end{axis}
	\end{tikzpicture}
	
			\caption{Impact of intercalated discs amplitude on CV for $k=3$.}
			\label{fig:plot_cv_vs_amp}
		\end{subfigure}
	\end{center}
	\caption{Effect of gap junctions surface area on CV.}
	\label{fig:intercalated_discs}
\end{figure}
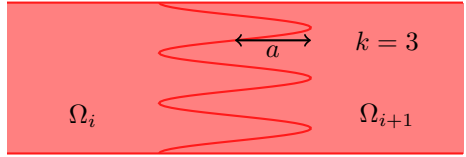
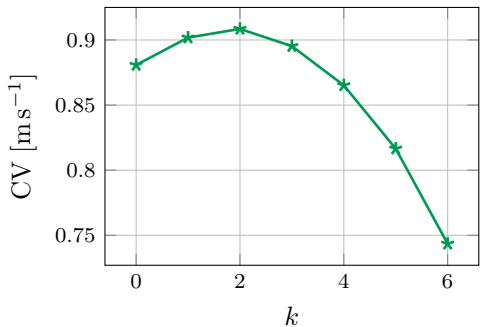
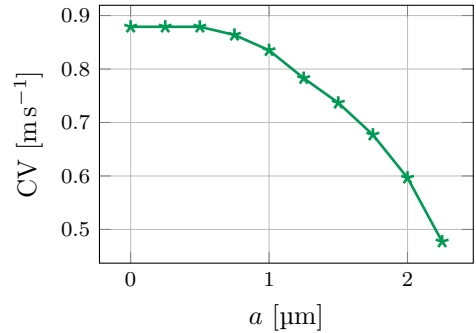

In general, gap junctions perpendicular to the fiber direction are not flat surfaces and are better modelled by intercalated discs \cite{HLR92}. In two dimensions we model these gap junctions with a sinus wave (see \cref{fig:ill_intercalated_discs} for an illustration) of amplitude $a$ and frequency $k$.

In this experiment we consider and array of $2\times 30$ cells of size $\SI{10}{\um}\times \SI{100}{\um}$ and measure CV as in \cref{sec:imp_dic_par_cv}. In \cref{fig:plot_cv_vs_freq} we show the conductive velocity as a function of the frequency $k$ for a fixed amplitude of $a=\SI{0.5}{\um}$. We note that for moderate frequency $k$ the CV increases due to an increase of contact surface area. However, for larger frequencies CV decreases, probably because of a flux saturation at the narrower junctions.
In \cref{fig:plot_cv_vs_amp} we show the conductive velocity as a function of the amplitude $a$ for a fixed frequency of $k=3$. Again, for larger amplitude $a$ the conduction velocity decreases.

\subsection{Dependence of conduction velocity on cells size and aspect ratio}\label{sec:CV_vs_cells_size}

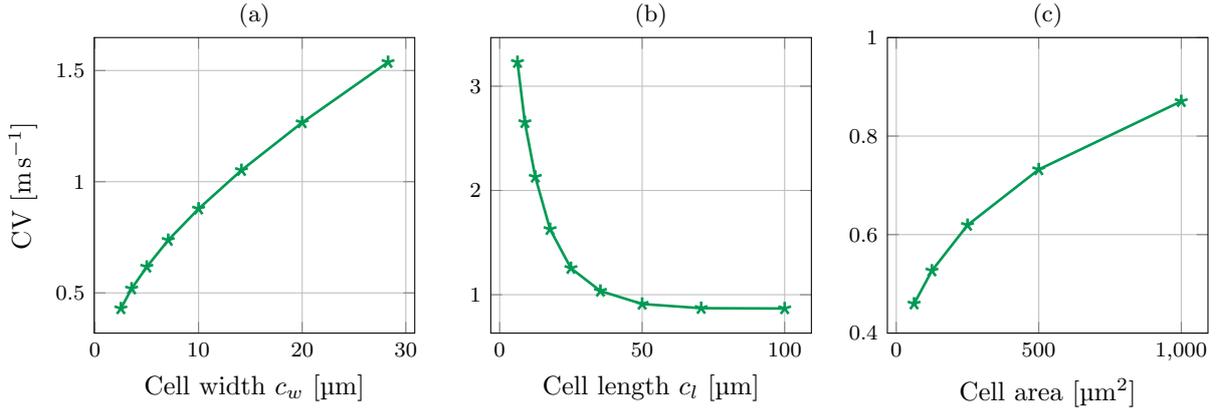
\begin{figure}[t]
\centering
\begin{tikzpicture}

\pgfplotscreateplotcyclelist{mycolorlisttest}{%
color=colorone,mark=\markone\\
color=colortwo,mark=\marktwo\\
color=colorthree,mark=\markthree\\
color=colorfour,mark=\markfour\\
color=colorfive,mark=\markfive\\
color=colorsix,mark=\marksix\\
}

\pgfplotsset{
    log x ticks with fixed point/.style={
    xticklabel={
      \pgfkeys{/pgf/fpu=true}
      \pgfmathparse{exp(\tick)}%
      \pgfmathprintnumber[fixed relative, precision=3]{\pgfmathresult}
      \pgfkeys{/pgf/fpu=false}
    }
  },
}

\begin{groupplot}[
    group style = {
        group size=3 by 1,
        horizontal sep=1cm,
        ylabels at = edge left,
    },
    grid = major,
    width = 5.8cm,
    height = 5.5cm,
    title style = {text width = 6cm},
    cycle list name = mycolorlisttest,
    legend cell align=left,
    legend pos = north west,
    title style={anchor=north, yshift=4ex},
    every axis plot/.append style=
    {
        solid,
        line width=\plotlinewidth pt,
        mark size=\plotmarksizeu pt,
    },
    tick label style = {font=\footnotesize},
    ylabel = {CV [\si{\meter\per\second}]},
]

    \nextgroupplot[
        title={\subcaption{\label{fig:cv_vs_cw}}},
        xlabel={Cell width $c_w$ [\um]},
    ]

    \addplot+ table [x=cw_um,y=CV,col sep=comma] {\currfiledir CV_vs_cw.csv};

    \nextgroupplot[
        title={\subcaption{\label{fig:cv_vs_cl}}},
        xlabel={Cell length $c_l$ [\um]},
    ]
    
    \addplot+ table [x=cl_um,y=CV,col sep=comma] {\currfiledir CV_vs_cl.csv};

    \nextgroupplot[
        title={\subcaption{\label{fig:cv_vs_clcw}}},
        xlabel={Cell area [\si{\micro\meter\squared}]},
        ymin = 0.4, ymax = 1,
    ]
    
    \addplot+ table [x=area_um,y=CV,col sep=comma] {\currfiledir CV_vs_clcw.csv};

\end{groupplot}
\end{tikzpicture}
\caption{Impact of cell length $c_l$, cell width $c_w$, and cell area with fixed aspect ratio on the CV.}
\label{fig:cv_vs_cellsize}
\end{figure}

Finally, we investigate how the cell's size and aspect ratio impact the conduction velocity. We consider an array of $2\times 30$ cells of size $c_w\times c_l$. First, we fix $c_w=\SI{10}{\um}$ and vary $c_l$, results are reported in \cref{fig:cv_vs_cl}. We observe as CV decreases as $c_l$ increases. In \cref{fig:cv_vs_cw} we display the results for fixed $c_l=\SI{100}{\um}$ and varying $c_w$, here CV increases with $c_w$. In the last figure \cref{fig:cv_vs_clcw} we vary both $c_l$ and $c_w$ while keeping a constant aspect ratio $c_l=10\cdot c_w$, more precisely they vary from $(c_w,c_l)=(\SI{2.5}{\um},\SI{25}{\um})$ to $(c_w,c_l)=(\SI{14}{\um},\SI{140}{\um})$. We see as CV increases with the cells area $A=c_l\cdot c_w$.

\section{Conclusion}\label{sec:conclusion}

In this paper we solve the cell-by-cell or EMI model for cardiac electrophysiology via the boundary element method, with no geometrical restrictions. The cell-by-cell model consists in Laplace equations inside and outside the cellular domains coupled with an ordinary differential equation on the transmembrane boundary and an algebraic condition on the gap junctions.
Due to the boundary integral formulation, Laplace equations are cleared away, yielding a differential algebraic equation living on the cell's boundaries only. In a subsequent step, the differential algebraic equation is reduced to an ordinary differential equation lying on the transmembrane boundary only.
Finally, we provide numerical results where: first, we study the accuracy of the numerical method and then we investigate the model properties and sensibility with respect to its parameters.

The convergence rate shows that the solution of the general problem is non-smooth, due to the presence of multiple cell contact (or 2 cells and the extracellular domain.) The single cell problem is instead smooth.
To the best of our knowledge, there are no regularity results for the single cell problem, except for those provided in \cite{matano2011} on asymptotic solutions. For the one-cell problem~\eqref{eq:cbcone} with smooth interface, but $\sigma_1 \neq \sigma_0$, the solution is probably regular. Intuitively, the interface problem with piecewise smooth coefficients and fixed transmembrane potential (that is, at the equilibrium) can be recast, via lifting~\cite{chen1998}, to a more classical interface problem already studied by \citet{babuska1970,kellogg1974}, who showed that the solution is at least $H^2$ on each subdomain. 

Well-posedness results for the general EMI problem are found in \cite{becue,franzone2002degenerate}, however we are not aware of regularity results for the general EMI problem~\eqref{eq:cbc}. For two or more cells in contact, subdomains must be polygonal, which limits the regularity. The singularities introduced by the contact have been analyzed by \citet{nicaise1994}. Moreover, the transmembrane voltage cannot be smooth on $\Gammam$, because it has multiple branches. These results should be taken in consideration in the development of higher order numerical schemes for the solution of the EMI problem.

We also show that the discretization parameters are not too restrictive, when compared to the more standard (homogenized) bidomain model. A typical time step is \SI{0.02}{\ms} or lower for IMEX solvers~\cite{KPD12}. In space, we observed here that a mesh resolution for the membrane of \SI{10}{\um} is sufficiently accurate for the cell-by-cell model. On the other hand, for the standard bidomain the mesh size depends on the front thickness, in turn depending on tissue excitability and conductivity. An accepted value is \SIrange{100}{200}{\um}~\cite{NKB11}. However, this is only true in the fiber direction and for healthy tissue, since in fibrotic tissue or in the cross-fiber direction the front thickness is generally lower~\cite{PHS16,fibernet}. Finally, we observed here that propagation failure can occur in the EMI model, in contrast to the bidomain model. This aspect is very important in the study of pathological situations.

This work paves the way in two directions. First, for designing another method where a more realistic cell-by-cell model is solved, i.e., where the linear algebraic condition on the gap junctions is replaced by a stiff nonlinear ordinary differential equation. Second, it provides the mathematical framework for solving the EMI model in three dimensions with a boundary integral formulation. Once assembled, the solution of the problem is very fast, opening interesting opportunities for long simulations.

\section*{Acknowledgement}
We are very grateful to Michael Multerer for kindly sharing the BEM code and the fruitful discussion on the problem. We also thank the MICROCARD consortium for the suggestions for improving the manuscript.

\bibliography{library}

\appendix 

\section{Poincaré--Steklov map on $\Omegae$}\label{app:PSe}
Here we briefly describe how to derive the Poincaré--Steklov operators $\PP_0,\PSe$ on $\Omega_0$ used in \cref{eq:contPSe,eq:PSe}, respectively.

Let $\bPhi\in\{\Gammam,\Sigma\}$, we introduce the restricted trace operators
\begin{equation}\label{eq:defgammaes}
	\begin{aligned}
		\gammae_{t,\bPhi} &\colon H^1(\Omegae) \to H^{1/2}(\bPhi), \quad &
		\gammae_{t,\bPhi} \ue(\bx) &= \lim_{\Omegae\ni \by\to \bx\in\bPhi} \ue(\by),  \\
		\gammae_{n,\bPhi} &\colon H^1(\Omegae) \to H^{-1/2}(\bPhi), &
		\gammae_{n,\bPhi} \ue(\bx) &= \lim_{\Omegae\ni \by\to \bx\in\bPhi}
		\langle \nabla \ue(\by), \bne \rangle.
	\end{aligned}
\end{equation}
From the Green's representation formula we have
\begin{equation}\label{eq:grfue}
	\begin{aligned}
		\ue(\bx) &
		= \int_{\Gammam} \gammae_{t,\by} G(\bx,\by) \gammae_{n,\Gammam} \ue(\by) \dif s_{\by}
		+\int_{\Sigma} \gammae_{t,\by} G(\bx,\by) \gammae_{n,\Sigma} \ue(\by) \dif s_{\by} &&\\
		&\quad - \int_{\Gammam} \gammae_{n,\by} G(\bx,\by) \gammae_{t,\Gammam} \ue(\by) \dif s_{\by}
		-\int_{\Sigma} \gammae_{n,\by} G(\bx,\by) \gammae_{t,\Sigma} \ue(\by) \dif s_{\by}, 
		&& \bx\in\Omegae.
	\end{aligned}
\end{equation}
For $\Phi,\Psi\in\{\Gammam,\Sigma\}$ , we define 
\begin{equation}
	\begin{aligned}
		\VVe^{\Psi,\Phi} &\colon H^{-1/2}(\Phi) \to H^{1/2}(\Psi), \quad &
		\VVe^{\Psi,\Phi} \rho(\bx) &= \int_{\Phi} \gammae_{0,\by}G(\bx,\by)\rho(\by)\dif s_{\by}, \quad & 
		\bx & \in\Psi, \\
		\KKe^{\Psi,\Phi} &\colon H^{1/2}(\Phi) \to H^{1/2}(\Psi), &
		\KKe^{\Psi,\Phi} \rho(\bx) &= \int_{\Phi} \gammae_{1,\by} G(\bx,\by) \rho(\by)\dif s_{\by}, & 
		\bx & \in\Psi, 
	\end{aligned}
\end{equation}
applying the trace operators $\gammae_{t,\Gammam}$ and $\gammae_{n,\Sigma}$ to \cref{eq:grfue} yields
\begin{equation}
	\begin{aligned}
		\gammae_{t,\Gammam}\ue &= 
		\VVe^{\Gammam,\Gammam}\gammae_{n,\Gammam}\ue
		+\VVe^{\Gammam,\Sigma}\gammae_{n,\Sigma}\ue
		-(\KKe^{\Gammam,\Gammam}-\tfrac{1}{2}I)\gammae_{t,\Gammam}\ue
		-\KKe^{\Gammam,\Sigma}\gammae_{t,\Sigma}\ue,\\
		\gammae_{t,\Sigma}\ue &= 
		\VVe^{\Sigma,\Gammam}\gammae_{n,\Gammam}\ue
		+\VVe^{\Sigma,\Sigma}\gammae_{n,\Sigma}\ue
		-\KKe^{\Sigma,\Gammam}\gammae_{t,\Gammam}\ue
		-(\KKe^{\Sigma,\Sigma}-\tfrac{1}{2}I)\gammae_{t,\Sigma}\ue,
	\end{aligned}
\end{equation}
which, after manipulation and setting $\gammae_{n,\Sigma}\ue=0$ (cf. \cref{eq:cbc_flux_ue_sigma}), result in
\begin{equation}\label{eq:gammae0Gammam}
	\begin{aligned}
		(\KKe^{\Gammam,\Gammam}+\tfrac{1}{2}I)\gammae_{t,\Gammam}\ue &= 
		\VVe^{\Gammam,\Gammam}\gammae_{n,\Gammam}\ue
		-\KKe^{\Gammam,\Sigma}\gammae_{t,\Sigma}\ue,\\
		(\KKe^{\Sigma,\Sigma}+\tfrac{1}{2}I)\gammae_{t,\Sigma}\ue &= 
		\VVe^{\Sigma,\Gammam}\gammae_{n,\Gammam}\ue
		-\KKe^{\Sigma,\Gammam}\gammae_{t,\Gammam}\ue.
	\end{aligned}
\end{equation}
Solving for $\gamma_{n,\Gammam}\ue,\gamma_{t,\Sigma}\ue$ with respect to $\gamma_{t,\Gammam}\ue$ yields the linear relation \cref{eq:contPSe} (dropping $\Gammam$ from the notation).

We discretize $\Gammam=\Gamma_1$ as in \cref{sec:cbem} (same collocation points) and place $\overline M$ collocation points $\overline\bx_j$ on $\Sigma$. We compute a smooth parametrization $\gamma_\Sigma:[0,1)\rightarrow \Rd$ of $\Sigma$, $\gamma_\Sigma(s_j)=\overline \bx_j$, by Fourier interpolation (as for $\Gammam$) and represent
\begin{equation}\label{eq:discgammasue}
	\begin{aligned}
		\gammae_{t,\Gammam}\ue (\gamma_{\Gammam}(t))&= \sum_{j=1}^{M} u_{0,\Gammam}^j L_j(t), & \gammae_{n,\Gammam}\ue (\gamma_{\Gammam}(t))&= \sum_{j=1}^{M} \tilde{u}_{0,\Gammam}^j L_j(t),\\
		\gammae_{t,\Sigma}\ue (\gamma_{\Sigma}(t))&= \sum_{j=1}^{\overline M} u_{0,\Sigma}^j \overline L_j(t),&
	\end{aligned}
\end{equation}
with $\overline L_j(s)$ the trigonometric Lagrange basis functions satisfying $\overline L_j(s_i)=\delta_{ij}$, $i,j=1,\ldots,\overline M$.
Inserting \cref{eq:discgammasue} into \cref{eq:gammae0Gammam} yields
\begin{equation}
	\begin{aligned}
		(\MKe^{\Gammam,\Gammam}+\tfrac{1}{2}I)\bm{u}_{0,\Gammam} &= 
		\MVe^{\Gammam,\Gammam}\bm{\tilde u}_{0,\Gammam}
		-\MKe^{\Gammam,\Sigma}\bm{u}_{0,\Sigma},\\
		(\MKe^{\Sigma,\Sigma}+\tfrac{1}{2}\overline I)\bm{u}_{0,\Sigma} &= 
		\MVe^{\Sigma,\Gammam}\bm{\tilde u}_{0,\Gammam}
		-\MKe^{\Sigma,\Gammam}\bm{u}_{0,\Gammam},
	\end{aligned}
\end{equation}
with $(\bm{u}_{0,\Gammam})_j=u_{0,\Gammam}^j,(\bm{\tilde u}_{0,\Gammam})_j=\tilde u_{0,\Gammam}^j,(\bm{u}_{0,\Sigma})_j=u_{0,\Sigma}^j$ and 
\begin{align}
	(\MKe^{\Gammam,\Gammam})_{kj} &= \KKe^{\Gammam,\Gammam}(L_j\circ\gamma_{\Gammam}^{-1})(\bx_k), & 
	(\MKe^{\Sigma,\Gammam})_{kj} &= \KKe^{\Sigma,\Gammam}(L_j\circ\gamma_{\Gammam}^{-1})(\overline\bx_k),\\
	(\MKe^{\Sigma,\Sigma})_{kj} &= \KKe^{\Sigma,\Sigma}(\overline L_j\circ\gamma_{\Sigma}^{-1})(\overline \bx_k), & 
	(\MKe^{\Gammam,\Sigma})_{kj} &= \KKe^{\Gammam,\Sigma}(\overline L_j\circ\gamma_{\Sigma}^{-1})(\bx_k),
\end{align}
and similarly for $\MVe^{\Sigma,\Gammam},\MVe^{\Gammam,\Gammam}$. Solving for $\bm{\tilde u}_{0,\Gammam},\bm{u}_{0,\Sigma}$ with respect to $\bm{u}_{0,\Gammam}$ yields 
\begin{equation}
	\bm{\tilde u}_{0,\Gammam}=\PSe \bm{u}_{0,\Gammam},
\end{equation}
which is employed in \cref{eq:PSe,eq:cbconedisc} (dropping $\Gammam$ from the notation).

\end{document}